\documentclass[english,ruled]{article}
\usepackage[T1]{fontenc}
\usepackage[latin9]{inputenc}
\usepackage{verbatim}
\usepackage{algorithm2e}
\usepackage{amsmath}
\usepackage{amsthm}
\usepackage{amssymb}
\usepackage{graphicx}
\usepackage{xcolor}
\usepackage{multicol}
\usepackage{multirow}
\usepackage{geometry}
\usepackage{enumitem}
\usepackage{setspace}
\makeatletter
\usepackage[toc,page,header]{appendix}
\usepackage{minitoc}
\usepackage{ifthen}
\newboolean{doublecolumn}
\setboolean{doublecolumn}{false}
\newboolean{arxiv}
\setboolean{arxiv}{true}
\usepackage{pgfplots}
\usepackage{pdflscape}
\pgfplotsset{compat=1.15}

\ifthenelse{\boolean{doublecolumn}}{
\renewcommand\vspace[1]{#1}
}

\renewcommand \thepart{}
\renewcommand \partname{}

\theoremstyle{plain}
\newtheorem{lem}{\protect\lemmaname}[section]
\theoremstyle{remark}
\newtheorem{rem}{\protect\remarkname}
\theoremstyle{plain}
\newtheorem{thm}{\protect\theoremname}[section]
\theoremstyle{plain}

\providecommand{\corollaryname}{Corollary}
\theoremstyle{plain}

\theoremstyle{plain}
\newtheorem{exple}{\protect\examplename}[section]

\usepackage{babel}
\providecommand{\propositionname}{Proposition}
\usepackage[colorlinks]{hyperref}       
\usepackage{url}            
\usepackage{booktabs}       
\usepackage{amsfonts}       
\usepackage{nicefrac}       
\usepackage{authblk}
\usepackage{optidef}

\newcommand\condspace{{\ifthenelse{\boolean{doublecolumn}}{}{\\}}}
\newcommand\condsep{{\ifthenelse{\boolean{doublecolumn}}{\vspace{-\topsep}}{}}}

\newcommand\tmop[1]{\mathrm{#1}}
\newcommand\tmtextbf[1]{\textbf{#1}}
\newcommand{\tmmathbf}[1]{\ensuremath{\boldsymbol{#1}}}

\makeatother
\usepackage{babel}

\usepackage[numbers]{natbib}
\usepackage[resetlabels]{multibib}
\usepackage[numbers]{natbib}

\newcites{app}{References in appendix}

\providecommand{\lemmaname}{Lemma}
\providecommand{\remarkname}{Remark}
\providecommand{\theoremname}{Theorem}
\providecommand{\examplename}{Example}

\newcommand{\assign}{:=}
\newcommand{\cdummy}{\cdot}

\newcommand{\tmstrong}[1]{\textbf{#1}}

\newcommand{\x}{\mathbf{x}}
\newcommand{\y}{\mathbf{y}}
\newcommand{\s}{\mathbf{s}}
\newcommand{\A}{\mathbf{A}}
\newcommand{\tmb}{\mathbf{b}}
\newcommand{\tmc}{\mathbf{c}}
\newcommand{\tmu}{\mathbf{u}}

\newcommand{\0}{\mathbf{0}}
\newcommand{\D}{\mathbf{D}}
\newcommand{\tmr}{\mathbf{r}}
\newcommand{\X}{\mathbf{X}}
\newcommand{\bs}{\mathbf{S}}
\newcommand{\tma}{\mathbf{a}}
\newcommand{\tmd}{\mathbf{d}}
\newcommand{\M}{\mathbf{M}}
\newcommand{\e}{\mathbf{e}}
\newcommand{\I}{\mathbf{I}}
\newcommand{\maxf}[1]{\underset{#1}{\text{maximize}}}
\newcommand{\minf}[1]{\underset{#1}{\text{minimize}}}

\newcommand{\w}{\mathbf{w}}
\newcommand{\z}{\mathbf{z}}
\newcommand{\tmv}{\mathbf{v}}
\newcommand{\U}{\mathbf{U}}

\newcommand{\W}{\mathbf{W}}
\newcommand{\Y}{\mathbf{Y}}
\newcommand{\Z}{\mathbf{Z}}
\newcommand{\B}{\mathbf{B}}
\newcommand{\tmP}{\mathbf{P}}

\begin{document}

\global\long\def\inprod#1#2{\left\langle #1,#2\right\rangle }%

\global\long\def\inner#1#2{\langle#1,#2\rangle}%

\global\long\def\binner#1#2{\big\langle#1,#2\big\rangle}%

\global\long\def\Binner#1#2{\Big\langle#1,#2\Big\rangle}%

\global\long\def\norm#1{\lVert#1\rVert}%

\global\long\def\bnorm#1{\big\Vert#1\big\Vert}%

\global\long\def\Bnorm#1{\Big\Vert#1\Big\Vert}%

\global\long\def\setnorm#1{\Vert#1\Vert_{-}}%

\global\long\def\bsetnorm#1{\big\Vert#1\big\Vert_{-}}%

\global\long\def\Bsetnorm#1{\Big\Vert#1\Big\Vert_{-}}%


\global\long\def\brbra#1{\big(#1\big)}%
\global\long\def\Brbra#1{\Big(#1\Big)}%
\global\long\def\rbra#1{(#1)}%


\global\long\def\sbra#1{[#1]}%
\global\long\def\bsbra#1{\big[#1\big]}%
\global\long\def\Bsbra#1{\Big[#1\Big]}%


\global\long\def\cbra#1{\{#1\}}%
\global\long\def\bcbra#1{\big\{#1\big\}}%
\global\long\def\Bcbra#1{\Big\{#1\Big\}}%
\global\long\def\vertiii#1{\left\vert \kern-0.25ex  \left\vert \kern-0.25ex  \left\vert #1\right\vert \kern-0.25ex  \right\vert \kern-0.25ex  \right\vert }%

\global\long\def\matr#1{\bm{#1}}%

\global\long\def\til#1{\tilde{#1}}%

\global\long\def\wtil#1{\widetilde{#1}}%

\global\long\def\wh#1{\widehat{#1}}%

\global\long\def\mcal#1{\mathcal{#1}}%

\global\long\def\mbb#1{\mathbb{#1}}%

\global\long\def\mtt#1{\mathtt{#1}}%

\global\long\def\ttt#1{\texttt{#1}}%

\global\long\def\dtxt{\textrm{d}}%

\global\long\def\bignorm#1{\bigl\Vert#1\bigr\Vert}%

\global\long\def\Bignorm#1{\Bigl\Vert#1\Bigr\Vert}%

\global\long\def\rmn#1#2{\mathbb{R}^{#1\times#2}}%

\global\long\def\deri#1#2{\frac{d#1}{d#2}}%

\global\long\def\pderi#1#2{\frac{\partial#1}{\partial#2}}%

\global\long\def\limk{\lim_{k\rightarrow\infty}}%

\global\long\def\trans{\textrm{T}}%

\global\long\def\onebf{\mathbf{1}}%

\global\long\def\zerobf{\mathbf{0}}%

\global\long\def\zero{\bm{0}}%


\global\long\def\Euc{\mathrm{E}}%

\global\long\def\Expe{\mathbb{E}}%

\global\long\def\rank{\mathrm{rank}}%

\global\long\def\range{\mathrm{range}}%

\global\long\def\diam{\mathrm{diam}}%

\global\long\def\epi{\mathrm{epi} }%

\global\long\def\relint{\mathrm{relint} }%

\global\long\def\inte{\operatornamewithlimits{int}}%

\global\long\def\cov{\mathrm{Cov}}%

\global\long\def\argmin{\operatornamewithlimits{arg\,min}}%

\global\long\def\argmax{\operatornamewithlimits{arg\,max}}%

\global\long\def\tr{\operatornamewithlimits{tr}}%

\global\long\def\dis{\operatornamewithlimits{dist}}%

\global\long\def\sign{\operatornamewithlimits{sign}}%

\global\long\def\prob{\mathrm{Prob}}%

\global\long\def\spans{\textrm{span}}%

\global\long\def\st{\operatornamewithlimits{s.t.}}%

\global\long\def\dom{\mathrm{dom}}%

\global\long\def\prox{\mathrm{prox}}%

\global\long\def\for{\mathrm{for}}%

\global\long\def\diag{\mathrm{diag}}%

\global\long\def\and{\mathrm{and}}%

\global\long\def\st{\mathrm{s.t.}}%

\global\long\def\dist{\mathrm{dist}}%

\global\long\def\Var{\operatornamewithlimits{Var}}%

\global\long\def\raw{\rightarrow}%

\global\long\def\law{\leftarrow}%

\global\long\def\Raw{\Rightarrow}%

\global\long\def\Law{\Leftarrow}%

\global\long\def\vep{\varepsilon}%

\global\long\def\dom{\operatornamewithlimits{dom}}%

\global\long\def\tsum{{\textstyle {\sum}}}%

\global\long\def\Cbb{\mathbb{C}}%

\global\long\def\Ebb{\mathbb{E}}%

\global\long\def\Fbb{\mathbb{F}}%

\global\long\def\Nbb{\mathbb{N}}%

\global\long\def\Rbb{\mathbb{R}}%

\global\long\def\extR{\widebar{\mathbb{R}}}%

\global\long\def\Pbb{\mathbb{P}}%

\global\long\def\Mrm{\mathrm{M}}%

\global\long\def\Acal{\mathcal{A}}%

\global\long\def\Bcal{\mathcal{B}}%

\global\long\def\Ccal{\mathcal{C}}%

\global\long\def\Dcal{\mathcal{D}}%

\global\long\def\Ecal{\mathcal{E}}%

\global\long\def\Fcal{\mathcal{F}}%

\global\long\def\Gcal{\mathcal{G}}%

\global\long\def\Hcal{\mathcal{H}}%

\global\long\def\Ical{\mathcal{I}}%

\global\long\def\Kcal{\mathcal{K}}%

\global\long\def\Lcal{\mathcal{L}}%

\global\long\def\Mcal{\mathcal{M}}%

\global\long\def\Ncal{\mathcal{N}}%

\global\long\def\Ocal{\mathcal{O}}%

\global\long\def\Pcal{\mathcal{P}}%

\global\long\def\Scal{\mathcal{S}}%

\global\long\def\Tcal{\mathcal{T}}%

\global\long\def\Wcal{\mathcal{W}}%

\global\long\def\Xcal{\mathcal{X}}%

\global\long\def\Ycal{\mathcal{Y}}%

\global\long\def\Zcal{\mathcal{Z}}%

\global\long\def\i{i}%


\global\long\def\abf{\mathbf{a}}%

\global\long\def\bbf{\mathbf{b}}%

\global\long\def\cbf{\mathbf{c}}%

\global\long\def\dbf{\mathbf{d}}%

\global\long\def\fbf{\mathbf{f}}%

\global\long\def\lambf{\bm{\lambda}}%

\global\long\def\alphabf{\bm{\alpha}}%

\global\long\def\sigmabf{\bm{\sigma}}%

\global\long\def\thetabf{\bm{\theta}}%

\global\long\def\deltabf{\bm{\delta}}%

\global\long\def\sbf{\mathbf{s}}%

\global\long\def\lbf{\mathbf{l}}%

\global\long\def\ubf{\mathbf{u}}%

\global\long\def\vbf{\mathbf{v}}%

\global\long\def\wbf{\mathbf{w}}%

\global\long\def\xbf{\mathbf{x}}%

\global\long\def\ybf{\mathbf{y}}%

\global\long\def\zbf{\mathbf{z}}%

\global\long\def\Abf{\mathbf{A}}%

\global\long\def\Ubf{\mathbf{U}}%

\global\long\def\Pbf{\mathbf{P}}%

\global\long\def\Ibf{\mathbf{I}}%

\global\long\def\Ebf{\mathbf{E}}%

\global\long\def\Mbf{\mathbf{M}}%

\global\long\def\Qbf{\mathbf{Q}}%

\global\long\def\Lbf{\mathbf{L}}%

\global\long\def\Pbf{\mathbf{P}}%


\global\long\def\abm{\bm{a}}%

\global\long\def\bbm{\bm{b}}%

\global\long\def\cbm{\bm{c}}%

\global\long\def\dbm{\bm{d}}%

\global\long\def\ebm{\bm{e}}%

\global\long\def\fbm{\bm{f}}%

\global\long\def\gbm{\bm{g}}%

\global\long\def\hbm{\bm{h}}%

\global\long\def\pbm{\bm{p}}%

\global\long\def\qbm{\bm{q}}%

\global\long\def\rbm{\bm{r}}%

\global\long\def\sbm{\bm{s}}%

\global\long\def\tbm{\bm{t}}%

\global\long\def\ubm{\bm{u}}%

\global\long\def\vbm{\bm{v}}%

\global\long\def\wbm{\bm{w}}%

\global\long\def\xbm{\bm{x}}%

\global\long\def\ybm{\bm{y}}%

\global\long\def\zbm{\bm{z}}%

\global\long\def\Abm{\bm{A}}%

\global\long\def\Bbm{\bm{B}}%

\global\long\def\Cbm{\bm{C}}%

\global\long\def\Dbm{\bm{D}}%

\global\long\def\Ebm{\bm{E}}%

\global\long\def\Fbm{\bm{F}}%

\global\long\def\Gbm{\bm{G}}%

\global\long\def\Hbm{\bm{H}}%

\global\long\def\Ibm{\bm{I}}%

\global\long\def\Jbm{\bm{J}}%

\global\long\def\Lbm{\bm{L}}%

\global\long\def\Obm{\bm{O}}%

\global\long\def\Pbm{\bm{P}}%

\global\long\def\Qbm{\bm{Q}}%

\global\long\def\Rbm{\bm{R}}%

\global\long\def\Ubm{\bm{U}}%

\global\long\def\Vbm{\bm{V}}%

\global\long\def\Wbm{\bm{W}}%

\global\long\def\Xbm{\bm{X}}%

\global\long\def\Ybm{\bm{Y}}%

\global\long\def\Zbm{\bm{Z}}%

\global\long\def\lambm{\bm{\lambda}}%

\global\long\def\alphabm{\bm{\alpha}}%

\global\long\def\albm{\bm{\alpha}}%

\global\long\def\taubm{\bm{\tau}}%

\global\long\def\mubm{\bm{\mu}}%

\global\long\def\yrm{\mathrm{y}}%

\global\long\def\tmstrong#1{{\bf #1}}%
\global\long\def\tmem#1{#1}%
\global\long\def\tmop#1{#1}%

\global\long\def\tmverbatim#1{#1}%
\global\long\def\revised#1{{#1}}%
\global\long\def\changed#1{\textcolor{blue}{#1}}%

\newtheorem{definition}{Definition}
\newcommand{\TODO}[1]{\textbf{\textcolor{red}{TODO: #1}}}
\font\myfont=cmr12 at 16.1pt
\title{When Does Primal Interior Point Method Beat Primal-dual in Linear Optimization?}
\setlength{\parindent}{0pt}
 \author[1]{Wenzhi Gao\thanks{gwz@stanford.edu}}
 \author[2]{Huikang Liu\thanks{hkliu2014@gmail.com}}
 \author[1,3]{Yinyu Ye\thanks{yyye@stanford.edu}}
  \author[1,3]{Madeleine Udell\thanks{udell@stanford.edu}}
 \affil[1]{ICME, Stanford University}
 \affil[2]{Antai School of Economics and Management, Shanghai Jiaotong University}
 \affil[3]{Management Science and Engineering, Stanford University}

\maketitle

\begin{abstract}
  The primal-dual interior point method (IPM) is widely regarded as the
  most efficient IPM variant for linear optimization. In this paper, we demonstrate that the improved stability of the pure primal IPM can allow speedups relative to a primal-dual solver,
  particularly as the IPM approaches convergence. 
  The stability of the primal scaling matrix makes it possible to accelerate each primal IPM step using fast preconditioned iterative solvers for the normal equations. Crucially, we identify properties of the central path that make it possible to stabilize the normal equations. Experiments on benchmark datasets demonstrate the efficiency of primal IPM and showcase its potential for practical applications in linear optimization and beyond.
\end{abstract}

\section{Introduction} \label{sec:intro}

We consider the standard-form primal and dual linear programs (LP) \\
\vspace{-15pt}

\begin{minipage}{.5\textwidth}
\begin{align}
\textup{(P)}\qquad \minf{\x} &~\quad ~\langle \tmc, \x \rangle \nonumber \\
    \text{subject to} & \quad\A \x = \tmb \nonumber \\
                      & ~~\quad\x \geq \0 \nonumber
\end{align}
\end{minipage}%
\begin{minipage}{.5\textwidth}
\begin{align} 
\textup{(D)}\qquad \maxf{(\y, \s)} & \quad ~~~\langle \tmb, \y \rangle \nonumber \\
    \text{subject to} & \quad \A^\top \y + \s = \tmc \nonumber \\
                      & ~~\quad \s \geq \0 \nonumber
\end{align}
\end{minipage}

\vspace{8pt}
with primal and dual feasible region $\mathcal{F}_p \assign \{ \x : \A \x = \tmb,
\x \geq \0 \}$ and $\mathcal{F}_d \assign \{ ( \y, \s )
: \A^{\top} \y + \s = \tmc, \s \geq 0 \}$.\\

Since Karmarkar developed the first polynomial time interior point method {\cite{karmarkar1984new}}, 
a substantial body of literature has explored both the theoretical \cite{wright1997primal,ye2011interior,roos2005interior} and practical aspects \cite{andersen1996implementation,adler1989implementation,karmarkar1991computational,mehrotra1992implementation} of IPMs. Over the years, IPM has evolved into a mature technology applicable to both linear and nonlinear optimization problems \cite{luenberger1984linear}.\\

IPM was first introduced in its primal form, known as primal potential reduction \cite{karmarkar1984new}, which operates solely in the primal space. However, with the evolution of IPM, several variants have emerged \cite{roos2005interior} that demonstrate improved robustness and efficiency compared to the original primal methods. Among these, the primal-dual IPM has stood the test of time and is widely regarded as the most efficient variant for linear and conic problems \cite{gondzio2012interior}. One of the few exceptions is in semidefinite
programming (SDP), where the dual potential reduction algorithm
{\cite{benson2008algorithm,gao2022hdsdp}} outperforms primal-dual
methods in sparse SDPs. Despite its initial significance, primal IPM now receives less attention in linear optimization. It is, therefore, interesting to ask:

\begin{center}
\textit{Are there problems on which primal IPM outperforms primal-dual IPM?}
\end{center}

This paper provides a positive answer to the question by clarifying the spectral properties of the primal normal matrices: in particular, its stability near convergence.
We demonstrate that this enhanced stability can be exploited to speed up later iterations of an IPM solver using the primal IPM iteration. 

\paragraph{Contributions.}
\begin{itemize}[leftmargin=10pt]
	\item We show that the normal matrix of the primal IPM is stable: when primal IPM converges, the normal matrix converges and need not become degenerate. 
 This property reduces the cost of solving the normal equations using a preconditioned iterative method such as conjugate gradient. 
 As a by-product, we obtain an $\Ocal(n^{3.5} + n^{2.5} \log^2 (1/\varepsilon) )$ complexity result for primal non-degenerate LPs. 
	\item We develop a hybrid primal/primal-dual IPM solver, using primal IPM to speed up later iterations of a primal-dual IPM solver. Extensive numerical experiments verify the effectiveness of our approach for linear optimization. We also studied the convergence pattern of IPM and established the empirical conditions under which primal IPM is likely to outperform primal-dual IPM. 
\end{itemize}

\subsection{Related Work}

\paragraph{Iterative solver-based IPM. }
Solving the normal equations $( \A \D^2 \A^{\top} ) \Delta \y = \tmr$,
where $\D$ is a diagonal matrix,
often presents the most serious computational obstacle for large-scale LPs.
One popular idea is to avoid factorizations by instead solving the normal equations 
using preconditioned iterative methods, leading to so-called matrix-free IPMs 
that do not require entrywise access to the matrix {\cite{gondzio2012matrix}}. 
Practically, IPMs using iterative solvers have shown empirical success in many applications \cite{zanetti2023new,cipolla2023proximal,velazco2010note,schork2020implementation,wang2000adaptive}.
On the theory side, literature on these methods must analyze the convergence of IPMs
with inexact subproblem solves \cite{bellavia1998inexact, freund1999convergence, gondzio2012interior, bellavia2019inexact}, as iterative solvers generally cannot reach the same accuracy as factorization-based solvers. 
An important branch of research focuses on the design of preconditioners for the IPM normal matrices \cite{bergamaschi2021new,gondzio2022general,chu2024randomized,chowdhury2020faster,bergamaschi2004preconditioning,bergamaschi2007inexact}. 

\paragraph{Arithmetic complexity of IPM. } Improving the arithmetic complexity of IPM has been a longstanding goal in the IPM literature \cite{lee2014path, vavasis1996primal, dadush2020scaling}, and often builds on tools from theoretical computer science. State-of-the-art complexity results for LP have been improved to the same order of matrix multiplication \cite{cohen2021solving, van2020deterministic}. These results generally ensure the stability of the primal-dual scaled projection matrix through a carefully designed data structure, which reduces the cost of solving for the IPM search directions. This paper similarly relies on the stability of the normal matrix but exploits the natural stability of the primal normal matrix instead of requiring extra machinery. 

\subsection{Notations} Throughout the paper we use $\langle \cdummy, \cdummy
\rangle$ to denote Euclidean inner product and $\| \cdummy \|$ to denote
Euclidean norm or matrix operator norm. Boldface letters $\A, \tma$ denote matrices and vectors.  $\0$ and $\e$ denote all-zero and all-one vectors and $\I$
denotes the identity matrix of proper size. Given a vector $\x \in
\mathbb{R}^n, \diag ( \x )$ is an $n \times n$ matrix with
$\x$ on the diagonal. Following the convention of IPM literature, we use $\X,
\bs, \D$ to denote $\diag ( \x )$, $\diag ( \s
)$ and $\diag ( \tmd )$ respectively. We denote by $\mathcal{F}^0_p
\assign \{ \x : \A \x = \tmb, \x > \0 \}$ and $\mathcal{F}^0_d
\assign \{  ( \y, \s ) : \A^{\top} \y + \s = \tmc, \s > 0
\}$ the primal and dual relative interior, respectively. Given a
symmetric positive definite matrix $\M$, $\kappa ( \M ) := \|\M \| \|\M^{-1} \| =
\lambda_{\max} ( \M ) / \lambda_{\min} ( \M )$ denotes
its Euclidean condition number and $\M^{1 / 2}$ denotes its square-root. Given two symmetric matrices $\A, \mathbf{B}$,   we define the semidefinite order $\A\succeq\mathbf{B}$  if and only if $\A - \mathbf{B}$ is positive semidefinite. Given any
$\x \in \Rbb^n_{++}$, we define $\| \x \|_{- \infty} \assign \min_{1\leq j \leq n} \{ x_j \}$ and $\|\tma\|_\x := \|\X^{-1} \tma \|$. Given an index set $\mathcal{B}$, we use $\A_{\mathcal{B}}, \x_{\mathcal{B}}$ to denote slicing of columns of $\A$, or slicing of coordinates in $\x$ over indices $\mathcal{B}$. We define gradient projection operator for a full row-rank matrix $\A$ as $\tmP_{\A}:= \I - \A^\top (\A\A^\top)^{-1}\A$. Unlike its definition in standard linear algebra, the term \textit{normal matrix} in this paper will refer to the left-hand-side matrix of the normal equation.

\section{Interior Point Method and Scaling Matrix} \label{sec:ipmscal}

This section describes the background and notation for interior point methods and 
introduces the primal IPM.

\subsection{Scaling Matrix for Interior Point Methods }

Given a primal-dual feasible pair $( \x, \y, \s ) \in \mathcal{F}^0_p
\times \mathcal{F}^0_d$, an interior point search direction is
determined by the optimality condition with linearized complementarity condition $\X \bs \e = \mu \e$:
\begin{equation}
  \left(\begin{array}{ccc}
    \A &  & \\
    & \A^{\top} & \I\\
    \D_{\X} &  & \D_{\bs}
  \end{array}\right) \left(\begin{array}{c}
    \Delta \x\\
    \Delta \y\\
    \Delta \s
  \end{array}\right) = \left(\begin{array}{c}
    \0\\
    \0\\
    \mu \e - \X \bs \e
  \end{array}\right). \label{eqn:aug}
\end{equation}
The first two equations ensure primal-dual feasibility, while the last one linearizes $\X \bs \e = \mu \e$.
There are several different ways of linearization which correspond to different \textit{scaling matrices} \cite{todd1994scaling}:
\begin{itemize}[leftmargin=10pt]
  \item \textit{Primal scaling} linearizes $\s - \mu\X^{-1}\e = \0$: $(\s + \Delta \s) - \mu (\X^{-1}\e - \X^{-2}\Delta \x ) = \0$, resulting in the system 
  \begin{equation} \label{eqn:pscal}
  	\mu \X^{- 1} \Delta \x + \X \Delta \s = \mu \e - \X \bs \e .
  \end{equation}
  \item \textit{Primal-dual scaling} linearizes $\X \bs \e= \mu \e$, resulting in the system 
  \begin{equation} \label{eqn:pdscal}
  	\bs \Delta \x + \X \Delta \s = \mu \e - \X \bs \e.
  \end{equation}
  \item \textit{Dual scaling} linearizes $\x -\mu\bs^{-1}\e = \0$: $(\x + \Delta \x) - \mu (\bs^{-1}\e - \bs^{-2}\Delta \s ) = \0$, resulting in the system
  \begin{equation} \label{eqn:dscal}
  	\bs \Delta \x + \mu \bs^{- 1} \Delta \s = \mu \e - \X \bs \e.
  \end{equation}
\end{itemize}
All three ways coincide when $\x$ and $\s$ exactly satisfy $\X \bs \e = \mu \e$; otherwise, they define different linear systems for IPM search directions. IPM software factorizes either the augmented matrix  
\[ \left(\begin{array}{cc}
     \D_{\bs}^{- 1} \D_{\X} & \A^{\top}\\
     \A & 
   \end{array}\right) \]
(formed by by eliminating $\Delta \s$ from \eqref{eqn:aug}), 
or the normal matrix $\M = \A \D_{\X}^{- 1} \D_{\bs} \A^{\top}$ 
(formed by eliminating $\Delta \x$ as well).
The \emph{primal, primal-dual, and dual normal matrices} are defined respectively as
\begin{align*}
  \M_{\X} \assign  \A \X^2 \A^{\top}
  \qquad
  \M_{\X, \bs} \assign \A \X \bs^{- 1} \A^{\top}  
  \qquad \M_{\bs}  \assign \A \bs^{- 2} \A^{\top}.
\end{align*}
Most IPM software uses primal-dual scaling, as years of
computational experience shows that primal-dual scaling performs best for linear optimization {\cite{lustig1994interior,wright1997primal}}.
When the IPM converges, some elements of $\X \bs^{- 1}$ approach $+\infty$, 
resulting in ill-conditioning of the primal-dual normal matrix $\M_{\X, \bs}$. 
Although the solution to the normal equations remains bounded so long as
$\x$ and $\s$ are in the interior {\cite{vavasis1996primal}}, 
the primal-dual normal matrix $\M_{\X, \bs}$ itself is unstable along IPM iterations. 
In this respect, the primal normal matrix $\M_{\X}$ has better properties: convergence of IPM
implies stability of $\M_{\X}$. 
Moreover, assuming the convergence of interior point iterates $\{ \x_k \}$ to $\x^{\star}$, 
primal IPM solves a sequence of converging linear systems $\M_{\X_k} \rightarrow
\M_{\X^{\star}}$, which are not necessarily ill-conditioned
{\cite{guler1993degeneracy}}. 

\subsection{Primal Scaling Interior Point Method}

The first polynomial-time IPM {\cite{karmarkar1984new}} is a primal-only method based on projective transformation and reducing a potential function. 
Since then, many variants of IPM have emerged. 
The two best-known variants, using a primal
scaling matrix, are primal potential reduction {\cite{ye1991n}} 
and primal logarithmic barrier methods {\cite{roos1992polynomial}}. 
This paper uses the primal logarithmic barrier method {\cite{roos1992polynomial}} as the working
algorithm. Primal IPM converges with minimal assumptions on the problem data:

\begin{enumerate}[leftmargin=35pt,itemsep=2pt,label=\textbf{A\arabic*:},ref=\rm{\textbf{A\arabic*}}]
	\item  $\mathcal{F}^0_p$ and $\mathcal{F}^0_d$ are nonempty
and $\A$ has full row rank. \label{A1}
\end{enumerate}

The assumption \ref{A1}, well-known as the interior point condition, is pervasive in the IPM literature. Given $\x \in \mathcal{F}_p^0$, we can derive the primal IPM by perturbing the KKT system,
replacing complementary slackness $\X \bs \e = \0$ by 
$\X \bs \e = \mu \e $, or equivalently, $\mu \X^{- 1} \e = \s$, resulting in the perturbed KKT system
\begin{align} \tag{KKT}
  \A \x ={} \tmb, \qquad
  \A^{\top} \y + \s ={}  \tmc,  \qquad
  \s - \mu \X^{- 1} \e ={}  \0. \label{eqn:kkt}
\end{align}
Then we can expand the nonlinear equation to first order to obtain the system \eqref{eqn:aug} with primal scaling \eqref{eqn:pscal}
\begin{align} \label{eqn:pkktstep}
  \A \Delta \x ={} & \0 \nonumber\\
  \A^{\top} \Delta \y + \Delta \s ={} & \0 \\
  \Delta \s + \mu \X^{- 2} \Delta \x ={} & \mu \X^{- 1} \e - \s.
  \nonumber
\end{align}

While \eqref{eqn:kkt} aligns with the derivation of primal-dual IPM \cite{wright1997primal}, another  natural motivation for primal IPM is through the logarithmic barrier: 
we replace  $\x \geq 0$ by the log-barrier penalty
\[ \textstyle \varphi_{\mu} ( \x ) \assign \langle \tmc, \x \rangle - \mu \sum_{i = 1}^n
   \log x_i . \]
At each iteration, we minimize the first-order expansion of $\varphi_{\mu}$
over an ellipsoid $\mathcal{E}_{\beta}(\x) \assign
\{ \tmd : \| \tmd \|_\x \leq \beta \}$:
\begin{eqnarray} \label{eqn:pipm-trustregion}
  \Delta \x = \argmin_{\A \tmd = \0} & \varphi_{\mu} ( \x ) +
  \langle \nabla \varphi_{\mu} ( \x ), \tmd \rangle \quad 
  \text{subject to} \quad 
  \| \tmd \|_\x \leq \beta,
\end{eqnarray}

where the scaled distance $\|\tmd\|_{\x} = \|\X^{-1} \tmd\|$. 
The ellipsoid, known as Dikin's ellipsoid \cite{dikin1967iterative}, guarantees the feasibility of $\x + \tmd$ if $\beta \leq 1$.
With the same choice of $\mu$, \eqref{eqn:pkktstep} and \eqref{eqn:pipm-trustregion} yield the same projected Newton direction
\begin{equation}
  \Delta \x = -\X (\I - \X\A^\top (\A\X^2\A^\top)^{-1}\A\X) ( \tfrac{1}{\mu} \X \tmc - \e ) = -\X \tmP_{\A\X} ( \tfrac{1}{\mu} \X \tmc - \e )
  \label{eqn:projnt}
\end{equation}
up to a multiplicative constant. Primal IPM starts from some $\mu = \mu_1$, takes a projected Newton direction, and slightly reduces $\mu$ in every iteration (\textbf{Algorithm \ref{alg:pipm}}).

\begin{algorithm}[h]
{\textbf{input} Initial point $\x_1 \in \mathcal{F}_p^0$, barrier parameter $\mu_1 > 0$, barrier reduction parameter $\tau \in (0,
1)$}\\
\For{k =\rm{ 1, 2,...}}{

\tmtextbf{compute} $\Delta \x_k = - \X_k \tmP_{\A\X_k}  ( \tfrac{1}{\mu_k} \X_k \tmc - \e )$

\tmtextbf{choose}\quad~$\alpha_k \in (0, 1)$ such that $\x_k + \alpha_k
\Delta \x_k \in \mathcal{F}_p^0$\\
\tmtextbf{update}\quad$\x_{k + 1} = \x_k + \alpha_k \Delta
\x_{k}$\quad and\quad$\mu_{k + 1} = (1 - \tau) \mu_k$
}
\caption{Primal interior point method \label{alg:pipm}}
\end{algorithm}

With proper initialization of $\x_1$ and configuration of $\{\mu_k\}$, it takes \textbf{Algorithm \ref{alg:pipm}}
$\mathcal{O} ( \sqrt{n} \log (1 / \varepsilon) )$ iterations to
find an $\varepsilon$-optimal $\x_{\varepsilon}$ such that $\x_{\varepsilon}
\in \mathcal{F}_p$ and $\langle \tmc, \x_{\varepsilon} \rangle \leq
\langle \tmc, \x^{\star} \rangle + \varepsilon$ \cite{ye2011interior, roos1992polynomial} . Moreover, an
$\varepsilon$-optimal dual solution $( \y_{\varepsilon},
\s_{\varepsilon} ) \in \mathcal{F}_d$ such that $\langle \tmb, \y_{\varepsilon} \rangle \geq \langle \tmb, \y^{\star} \rangle - \varepsilon$ can be computed by a projection subproblem \cite{roos1992polynomial, ye1991n}.

\subsection{Primal Central Path}
To make the discussion self-contained, we review a few definitions concerning the primal central path {\cite{gonzaga2004properties}}. 

\begin{definition}[Primal central path \cite{roos1997theory}]
  The primal central path is defined by
  \[ \mathcal{C}_p \assign \{ \x_{\mu} : \mu > 0 \}, \qquad \x_{\mu} \assign \argmin_{\x \in \mathcal{F}_p} ~ \{
  \langle \tmc, \x \rangle - \mu \textstyle \sum_{i = 1}^n \log x_i \}.\]
\end{definition}
Under \ref{A1}, the primal central path is guaranteed to exist \cite{roos1997theory}.
Primal IPM (\tmtextbf{Algorithm \ref{alg:pipm}})
explicitly follows $\mathcal{C}_p$: 
the algorithm generates iterates that stay near $\mathcal{C}_p$
in the sense of the following proximity measure.

\begin{definition}[Proximity to central path  {\cite{roos1997theory}}] \label{def:proximity}
 Given $\x \in \mathcal{F}_p^0$, its centrality proximity measure $\delta ( \x, \mu )$ with respect to barrier
  parameter $\mu > 0$ is defined by
\begin{equation} \label{eqn:def-proximity}
	\delta ( \x, \mu ) := \min_{( \y, \s ) \in
     \mathcal{F}_d} \big\| \tfrac{1}{\mu} \X \s - \e \big\| = \|\tmP_{\A\X} ( \tfrac{1}{\mu} \X \tmc - \e ) \|,
\end{equation}
where $\tmP_{\A\X}$ is the same as in \eqref{eqn:projnt}.
\end{definition}
Evaluating central path proximity is essentially a projection problem in the
dual space, which minimizes the distance to $\mathcal{C}_p$. 
The proximity measure $\delta ( \x,
\mu )$ is a ``primal-only'' version of the $\ell_2$-neighborhood widely
used in primal-dual path-following IPM {\cite{wright1997primal}}. Finally,
we introduce the concept of maximal complementarity pair.

\begin{definition}[Maximal complementarity pair  {\cite{guler1993convergence}}] \label{def:maxcomp} 
Under \ref{A1}, there exists some
  partition $(\mathcal{B}, \mathcal{N})$ satisfying $\mathcal{B} \cup \mathcal{N}= \{ 1, \ldots, n \}$ and
  $\mathcal{B} \cap \mathcal{N}= \varnothing$ such that
  \begin{align}
    \mathcal{B}={} & \{ j : x^{\star}_j > 0 \text{~~for some~~}  \x^{\star} \},\qquad
    \mathcal{N}={}  \{ j : s^{\star}_j > 0 \text{~~for some~~}  \s^{\star}
    \}, \nonumber
  \end{align}
  
  where $( \x^{\star}, \s^{\star} )$ is some optimal
  primal-dual solution pair. 
\end{definition}
Partition $(\mathcal{B}, \mathcal{N})$ is known as
  the maximal complementarity pair.
  Moreover, $\x_{\mathcal{N}}^{\star} = \0$ and $\s_{\mathcal{B}}^{\star} = \0$ hold for
  any optimal solution $( \x^{\star}, \s^{\star} )$, and there exists some strictly
  complementary solution $( \x^{\star}, \s^{\star} )$ such that $\x^{\star} + \s^{\star} > \mathbf{0}$. We define constant
\begin{equation} \label{eqn:gamma}
  	\gamma :=  \tfrac{n}{\| \s^{\star}_{\mathcal{N}} \|_{-
  \infty}} > 0
  \end{equation}
    for an arbitrarily chosen strictly complementary solution $( \x^{\star}, \s^{\star} )$.
\section{Interior Point Method in the Primal Form} \label{sec:pipm}

The bottleneck of primal IPM lies in computing the IPM direction \eqref{eqn:projnt}, where a normal equation with left-hand-side normal matrix $\M_{\X_k}$ has to be solved. As previously remarked, convergence of the IPM iterates $\{\x_k\}$ implies the convergence of $\{\M_{\X_k}\}$. In this section, we show it is possible to reduce the cost of solving normal equations by exploiting the fact that solving a sequence of linear systems is easier.

\subsection{Convergence of the Scaling Matrix}

We start by characterizing the convergence behavior of primal IPM under primal non-degeneracy:

\begin{enumerate}[leftmargin=35pt,itemsep=2pt,label=\textbf{A\arabic*:},ref=\rm{\textbf{A\arabic*}},start=2]
	\item The primal problem is non-degenerate; the maximal complementarity partition satisfies $| \mathcal{B} | = m$. \label{A2}
\end{enumerate}

The primal optimal solution $\x^\star$ is unique given \ref{A2}. The proof is simple: we show that for a fixed proximity parameter
$\lambda$, the diameter of the proximity set $\mathcal{P}_{\mu} \assign \{ \x : \delta
( \x, \mu ) \leq \lambda \}$ is bounded by $\mathcal{O}
(\mu)$. Since primal IPM decreases $\mu$ linearly and its iterates stay within the proximity set, we obtain
convergence of $\{ \x_k \}$ as $\mu \rightarrow 0$.

\begin{lem} \label{lem:boundedprox}
Under \ref{A1} and \ref{A2}, for any $\lambda \in [ 0, \tfrac{\sqrt{2}}{2} ]$ and $\mu > 0$, the 
  central path proximity set
 $\mathcal{P}_{\mu} = \{ \x : \delta ( \x, \mu )
     \leq \lambda \}$,
satisfies ${\diam} (\mathcal{P}_{\mu}) =\mathcal{O} (\mu)$.
\end{lem}

With \textbf{Lemma \ref{lem:boundedprox}}, we obtain a global iteration convergence rate for primal IPM.

\begin{thm} \label{thm:iterconv}
  Under \ref{A1} and \ref{A2}, suppose an algorithm
  maintains $\delta ( \x_k, \mu_k ) \leq \frac{\sqrt{2}}{2}$ for all $k = 1, 2, \ldots$ and that
  $\mu_{k + 1} = (1 - \tau) \mu_k$ for $\tau \in (0, 1/\sqrt{2 n}]$. Then $\| \x_i - \x_j \| =\mathcal{O} (\mu_k)$ for all $i, j \geq k$.
\end{thm}

Since the primal normal matrix $\M_{\X} = \A \X^2 \A^{\top}$ consists of quadratic terms with respect to
$\x$, convergence of $\{ \x_k \}$ implies
convergence of $\{ \M_{\X_k} \}$. Therefore, the final iterations of the primal IPM solve linear systems with converging coefficient matrices: normal matrices of primal IPM stay in proximity to each other. This observation opens a way for preconditioned iterative solvers in primal IPM.

\subsection{Proximity in Euclidean Distance}
 Intuitively, two nearby non-degenerate matrices have similar spectra. \textbf{Lemma \ref{lem:euc}} quantifies this intuition by analyzing the conditioning of the preconditioned systems.
 
 \begin{lem}[Conditioning within the Euclidean ball] \label{lem:euc}
  Let $\lambda_{\x} \assign \lambda_{\min} ( \M_{\X} ) > 0$. 
  For any $\x \in \mathbb{R}_{+ +}^n$, $\beta \leq 1$ and $\Delta \x$ such
  that $\| \Delta \x \| \leq
  \frac{\sqrt{\lambda_{\x}}}{\| \A \|} \beta$,  we have $(1 - \beta)^2 \cdot \I \preceq \M_{\X}^{- 1 / 2} ( \M_{\X + \Delta \X} ) \M_{\X}^{- 1 / 2} \preceq (1 + \beta)^2 \cdot \I $ and 
  \[ \kappa ( \M_{\X}^{- 1 / 2} ( \M_{\X + \Delta \X} )
     \M_{\X}^{- 1 / 2} ) \leq \big( \tfrac{\beta + 1}{\beta - 1}
     \big)^2. \]
\end{lem}

\textbf{Lemma \ref{lem:euc}} shows that when $\M_\X$ is far from singular and $\Delta \x$ is small, we can solve $\M_{\X + \Delta\X}$ using $\M_{\X}$ as a preconditioner: proximity in Euclidean distance implies effective preconditioning. In particular, \textbf{Lemma \ref{lem:euc}} and \textbf{Theorem \ref{thm:iterconv}} together imply that when the primal IPM converges, good conditioning is guaranteed for all the normal matrices generated by the primal IPM iteration:
\begin{thm} \label{thm:condnum}
  Under \ref{A1} and \ref{A2}, let $\{ \M_{\X_{k}} \}$ be
  the sequence of normal matrices generated by the primal IPM specified in \textbf{Algorithm \ref{alg:pipm}} and suppose the conditions in \textbf{Theorem \ref{thm:iterconv}} hold.
Then
  \[ \kappa ( \M_{\X_{k}}^{- 1 / 2} \M_{\X_j} \M_{\X_{k}}^{- 1 / 2} )
     \leq 9 \]
  for all $j \geq k \geq k^\star \assign \min \big\{ k: \mu_{k} \leq \tfrac{\sqrt{\lambda^{\star}} \tau}{16
     \gamma \|\A\|} \cdot ( 1 + \tfrac{\|
     \A_{\mathcal{N}} \|}{\sqrt{\lambda_{\min} (
     \A_{\mathcal{B}}^{\top} \A_{\mathcal{B}} )}})^{-1} \big\}$, where $\lambda^\star := \lambda_{\min} ( \M_{\X^\star} )$. In particular, choosing $\tau = 1/\sqrt{2n}$ gives \[k^\star = \Ocal\Big(\sqrt{n} \Big[\log (\tfrac{\gamma \|\A\|}{\sqrt{\lambda^\star}}) + \log\big(1+\tfrac{\|\A_{\Ncal}\|}{\sqrt{\lambda_{\min} (\A_\Bcal^\top \A_\Bcal)}}\big)\Big] \Big).\]
     
\end{thm}
Given a linear system with a bounded condition number after preconditioning, it is well-known that preconditioned iterative solvers, such as conjugate gradient, can solve the linear system to $\varepsilon$ accuracy in $\Ocal(n^2 \log(1/\varepsilon))$ arithmetic complexity \cite{saad2003iterative}. \textbf{Theorem \ref{thm:condnum}} itself provides an algorithm: when the primal IPM iterates converge, we reuse the factorization of past normal matrices as preconditioners to solve future normal equations. 

\begin{algorithm}[h]
{\textbf{input} Initial points $\z=\x_1 \in \mathcal{F}_p^0$,  barrier parameter $\mu_1 > 0$, barrier reduction parameter $\tau \in (0,
1),$}\\ \qquad~~~ {cached factorization of $\M_\Z = \A\Z ^2\A^\top,$ factorization update tolerance $\theta > 0$}\\
\For{k =\rm{ 1, 2,...}}{

\If{$\|\x_k - \z\| \geq \theta$}{
\tmtextbf{update} $\z = \x_k$ \textbf{and factorize}  $\M_\Z = \A\X_k ^2\A^\top$\\
}
\tmtextbf{compute} $\Delta \hat{\x}_{k}$ with  iterative solver and preconditioner $\M_\Z$

\tmtextbf{choose}\quad~$\alpha_k \in (0, 1]$ such that $\x_k + \alpha_k
\Delta \hat{\x}_{k} \in \mathcal{F}_p^0$\\
\tmtextbf{update}\quad$\x_{k + 1} = \x_k + \alpha_k \Delta
\hat{\x}_{k}$\quad and\quad$\mu_{k + 1} = (1 - \tau) \mu_k$
}
\caption{Preconditioned primal interior point method \label{alg:pipm-precond}}
\end{algorithm}

Combining \textbf{Theorem \ref{thm:condnum}} with the IPM inexact analysis yeilds the arithmetic complexity for \textbf{Algorithm \ref{alg:pipm-precond}}.

\begin{thm}[Informal] \label{thm:informalinexact}
   The total complexity of \tmtextbf{Algorithm \ref{alg:pipm-precond}} is $\mathcal{O} ( n^{3} k^{\star} + n^{2.5} \log^2 (
     \tfrac{1}{\varepsilon} ) )$, where $k^\star$ is defined in \textbf{Theorem \ref{thm:condnum}}.
\end{thm}

 The proof of \tmtextbf{Theorem \ref{thm:informalinexact}} reveals the 
 algorithm has two distinct phases. The $n^{3} k^{\star}$ complexity follows from the standard analysis of primal IPM \cite{ye2011interior}.
When $\mu_k$ is small enough that all future normal matrices $\{ \M_k \}_{k \geq k^{\star}}$ are close to $\M_{k^{\star}}$, 
using the fixed preconditioner $\M_{k^{\star}}$ decreases
 the cost per IPM iteration from $\mathcal{O} (n^3)$ to $\mathcal{O}(n^2)$. 
  Proof of \tmtextbf{Theorem \ref{thm:informalinexact}} is deferred to \textbf{Appendix \ref{app:inexact}}.
  
\tmtextbf{Algorithm \ref{alg:pipm-precond}} is simple and straight-forward to implement. Yet in practice, $\mu_{k^{\star}}$ can be extremely small, so the algorithm may terminate before the second, faster, phase of convergence begins. To make primal IPM practically useful, we must introduce one more type of proximity defined with respect to the scaled distance.

\subsection{Proximity in the Scaled Distance}
Given $\x, \w \in \Rbb^n_{++}$, we say $\w$ is close to $\x$ in the scaled distance if $\|\w - \x\|_\x$ is small. An important consequence of proximity in the scaled distance is \textbf{Lemma \ref{lem:shiftedscal-l2}}, known as the \textit{shifted scaling theorem} in the IPM literature. It shows that given $\w = \x + \Delta \x$ that is close to $\x$ in the scaled distance, the scaled projection operator $\tmP_{\A\X}$ can be well approximated by $\tmP_{\A\W}$.

\begin{lem}[Shifted scaling {\cite{gonzaga2004properties}}] \label{lem:shiftedscal-l2} For any $\x \in \Rbb^n_{++}, \beta \leq \frac{1}{4}$ and $\Delta \x$ such that $\| \Delta \x \|_{\x} \leq \beta$, 
  \begin{align}
    & \big\| \big[ \X^{- 1}(\X +\Delta \X) \tmP_{\A(\X + \Delta \X)}  \X^{- 1}(\X +\Delta \X) - \tmP_{\A\X} \big]
    \tmv \big\|
    \leq{}  3 \beta \| \tmP_{\A\X} \tmv \| ~~  \text{for all} ~~ \tmv \in
  \mathbb{R}^n. \nonumber
  \end{align}
\end{lem}

\textbf{Lemma \ref{lem:shiftedscal-l2}} was initially proposed as a theoretical tool to analyze the primal central path \cite{gonzaga2004properties}.  However, it will be a critical tool for our algorithm design: suppose we are at some $\x \in \Fcal_p^0$. Then, according to \eqref{eqn:projnt}, the primal IPM direction at $\x$ is \[ -\X \tmP_{\A\X} ( \tfrac{1}{\mu} \X \tmc - \e ).\]

Suppose instead of applying $\tmP_{\A\X} = \I - \X\A^\top (\A\X^2\A^\top)^{-1}\A\X$, which involves $(\A\X^2\A^\top)^{-1}$, we only have access to the operator $\tmP_{\A\W}$ at another point $\w$ such that $\|\w - \x\|_\x$ is small. Then, using \textbf{Lemma \ref{lem:shiftedscal-l2}} with $\x + \Delta \x = \w$ and $\tmv = \tfrac{1}{\mu} \X \tmc - \e$,  we notice that
\begin{align}
{} & \big\| \big[-\X \tmP_{\A\X} (\tfrac{1}{\mu} \X \tmc - \e) \big] -\big[ - \W \tmP_{\A\W}  \X^{- 1}\W (\tfrac{1}{\mu} \X \tmc - \e) \big]
     \big\|_\x \nonumber \\
	= {}& \big\|  \X^{- 1}\W \tmP_{\A\W}  \X^{- 1}\W (\tfrac{1}{\mu} \X \tmc - \e) - \tmP_{\A\X} (\tfrac{1}{\mu} \X \tmc - \e)
     \big\| \label{eqn:scal-dist}  \\
   ={} & \big\| \big[ \X^{- 1}\W \tmP_{\A\W}  \X^{- 1}\W - \tmP_{\A\X} \big]
    (\tfrac{1}{\mu} \X \tmc - \e) \big\|
    \leq{} 3 \beta  \| \tmP_{\A\X} (\tfrac{1}{\mu} \X \tmc - \e) \|= 3 \beta \delta(\x, \mu), \label{eqn:shift-scal-prox}
\end{align}
where \eqref{eqn:scal-dist} uses the definition $\|\tma\|_\x = \|\X^{-1} \tma \|$ and \eqref{eqn:shift-scal-prox} uses $\delta(\x, \mu) = \| \tmP_{\A\X} (\tfrac{1}{\mu} \X \tmc - \e) \|$ from \eqref{eqn:def-proximity}. Hence, if $\delta(\x, \mu)$ is small and $\w$ is close to $\x$ in the scaled distance, the shifted scaling theorem provides a good surrogate direction 
\begin{equation} \label{eqn:surrodirection}
- \W \tmP_{\A\W}  \X^{- 1}\W (\tfrac{1}{\mu} \X \tmc - \e)	
\end{equation}
 without even computing $\tmP_{\A\X}$. This observation is significant in making primal IPM practically useful.\\

In the next section, we put all the things together and
show how these two types of proximity can be combined to make primal IPM competitive in solving real LPs.

\subsection{Combining two Types of Proximity}

We have shown that proximity in Euclidean and scaled distance both allow us to reduce the cost of solving the normal equations. However, as the following example shows, neither of them \textit{alone} may be usable in practice.

\begin{exple} \label{exp:1}
Let $\z = (10^{10}, 10^{- 5})^\top$ and $\x = (10^{10} - 10^5, 10^{-
10})^\top$. Then
\[  \quad \| \z - \x 
   \| \approx 10^5 \quad \text{and} \quad \|   \z - \x
   \|_\x \approx 10^5 . \]	
 $\x$ and $\z$ are far away in both Euclidean and the scaled distance. 
\end{exple}

\textbf{Example \ref{exp:1}} can indeed happen in a real algorithm: the first coordinate corresponds to a basic variable converging to a large value, and the second is a non-basic variable converging to 0. Although $\x$ and $\z$ are both near convergence, neither Euclidean nor the scaled distance is small. The issue arises due to the misalignment of two types of distances: as the primal iterates converge, the coordinates converging to 0 will be close in Euclidean distance but not the scaled distance. In contrast, the coordinates converging to a (large) nonzero value will be close in the scaled distance long before they are close in Euclidean distance. Two geometries complement each other. Hence, we define a combined metric to reap the benefits of both.

\begin{definition}[Thresholded scaled distance]
  The $\nu$-thresholded scaled distance between $\y, \z \in \mathbb{R}^n$
  with respect to $\x \in \mathbb{R}^n_{++}$ is given by
  \[ \| \y - \z \|_{\x, \nu} = \sqrt{\| \X_{\mathcal{L}}^{-
     1} ( \y_{\mathcal{L}} - \z_{\mathcal{L}} ) \|^2 +
     \| \y_{\mathcal{S}} - \z_{\mathcal{S}} \|^2}, \]
  where $\mathcal{L} \assign \{ j : x_j \geq \nu \}$ and $\mathcal{S}= \{ j
  : x_j < \nu \}$.
\end{definition}

The distance $\| \y - \z \|_{\x, \nu}$ combines Euclidean and the scaled
distance with a thresholding parameter $\nu$: distance with respect to
large elements $\mathcal{L}$ are scaled down, while the small elements $\mathcal{S}$ are not. For example, in \textbf{Example \ref{exp:1}}, taking $\nu = 1$ gives $\|\z - \x \|_{\x, 1} \approx 10^{-5}$. This new metric motivates a modification of primal IPM that takes advantage of both geometries.

\subsubsection{Primal Interior Point Method with Delayed Scaling}

We are ready to leverage the thresholded scaled distance to motivate a more practical primal IPM.  Given $\x, \z \in \Fcal_p^0$ 
and a factorization for $\M_{\Z}$, 
suppose we want to
compute the primal search direction
\[ \Delta \x = - \X \tmP_{\A\X} ( \tfrac{1}{\mu} \X \tmc - \e ) . \]
Following our previous argument, we hope the factorization of $\M_{\Z}$ can accelerate computing $\Delta \x$. But even if $\| \z - \x
\|_{\x, \nu}$ is small, there is no guarantee that either $\| \z - \x \|$ or $\| \z - \x \|_\x$ is small.
Instead, our algorithm (\textbf{Algorithm \ref{alg:delayed}}) constructs an intermediate point $\w$ so that both $\| \w -
\z \|$ and $\|\w - \x\|_\x$ are small. Then we compute a surrogate search direction $\Delta \hat{\x} = - \W \tmP_{\A\W} ( \tfrac{1}{\mu} \W \tmc - \W\X^{-1}\e )$ in \eqref{eqn:surrodirection} with preconditioner $\M_\Z$.

\begin{algorithm}[h]
{\textbf{input} Initial point $\z = \x_1 \in \mathcal{F}_p^0$, barrier parameter $\mu_1 > 0$, barrier reduction parameter  $\tau \in (0,
1)$,\\
\qquad~~~ cached factorization of $\M_\Z = \A\Z ^2\A^\top$, threshold  $\nu > 0$, factorization update tolerance $\theta > 0$}\\
\For{k =\rm{ 1, 2,...}}{

\If{$\|\x_k - \z\|_{\x_k, \nu} \geq \theta $}{
\tmtextbf{update} $\z = \x_k$ \tmtextbf{and factorize} $\M_\Z = \A\Z ^2\A^\top$\\
}
\tmtextbf{choose} delayed scaling matrix $\w$, $\w_\mathcal{S} = \x_{k,\Scal}, \w_{\Lcal} = \z_{\Lcal}$ 

\tmtextbf{compute} $\Delta \hat{\x}_{k} = - \W \tmP_{\A\W}  ( \tfrac{1}{\mu_k} \W \tmc - \W\X^{-1}_k\e )$ with preconditioner $\M_\Z$

\tmtextbf{choose}\quad~$\alpha_k \in (0, 1)$ such that $\x_k + \alpha_k
\Delta \hat{\x}_{k} \in \mathcal{F}_p^0$\\
\tmtextbf{update}\quad  $\x_{k + 1} = \x_k + \alpha_k \Delta
\hat{\x}_{k}$\quad and\quad$\mu_{k + 1} = (1 - \tau) \mu_k$
}
\caption{Preconditioned primal interior point method with delayed scaling matrix \label{alg:delayed}}
\end{algorithm}

The validity of \textbf{Algorithm \ref{alg:delayed}} is guaranteed by \textbf{Theorem \ref{thm:delayedscal}}.
\begin{thm}\label{thm:delayedscal} 
Fix $\nu \in (0, 1]$. Suppose $\mu > 0$, $\x, \z \in \mathcal{F}^0_p$ and $\M_{\Z}
  \succeq \lambda_{\z} \cdummy \I$ is already factorized and that 
  \[ \| \z - \x \|_{\x, \nu} \leq \min \big\{ \tfrac{\sqrt{\lambda_{\z}}}{2 \| \A \|}, \tfrac{1}{4}, \tfrac{
     \varepsilon}{6 \delta ( \x, \mu )} \big\}.  \]
Then \textbf{Algorithm \ref{alg:delayed}} computes $\Delta \hat{\x}$ such that $\| \Delta \x - \Delta \hat{\x}  \|_\x \leq
  \varepsilon$ in $\mathcal{O} (n^2 \log (1 / \varepsilon))$ arithmetic complexity with the cached factorization of $\M_{\Z}$.
\end{thm}

To illustrate the choice of $\w$, again in \textbf{Example \ref{exp:1}}, we have

\begin{equation*}
	\x = \footnotesize\begin{pmatrix}
     ~10^{10}-10^5~\\
     10^{-10}
   \end{pmatrix} \overset{\text{Scaled}}{\approx} \w =
 \begin{pmatrix}
     ~10^{10}~\\
     10^{- 10}
   \end{pmatrix} \overset{\text{Euclidean}}{\approx} \begin{pmatrix}
     ~10^{10}~\\
     10^{-5}
   \end{pmatrix} = \z 
\end{equation*}\\
with $\| \w - \z \| \approx 10^{- 5}, \|\w -
\x \|_\x \approx 10^{- 5}$. 
More generally, good preconditioning requires some scaling matrix $\W$ such that \textbf{1)} $\w$ is close to $\z$ in Euclidean distance so that \textbf{Lemma \ref{lem:euc}} guarantees good preconditioning. \textbf{2)} $\w$ is close to $\x$ in the scaled distance so that \textbf{Lemma \ref{lem:shiftedscal-l2}} guarantees that the search direction is still accurate. Such a $\w$ exists in primal IPM, and we call it a \emph{delayed} scaling matrix: the small coordinates in $\x$ are updated to form $\w$, but the update of the large coordinates in $\x$ is delayed to maximize the stability of the normal matrices, by \textbf{Lemma \ref{lem:euc}}.

\section{Primal Interior Point Method in Practice}
\label{sec:pipmprac}

This section describes how to implement the primal
IPM, including an infeasible start extension and dualization,
and shows how to generalize primal IPM to SDP.

\subsection{Infeasible Start Primal IPM}

Most modern IPM solvers use an infeasible start variant of primal-dual IPM. Infeasible start IPM only requires iterates to lie in $\Rbb_{++}^n$, while
primal and dual feasibility can be violated. We derive an infeasible start
variant of primal IPM through the KKT interpretation \eqref{eqn:pkktstep} in \tmtextbf{Section \ref{sec:ipmscal}}.\ Given a primal-dual solution pair $( \x, \y, \s )$, primal, dual
and complementary residuals are, respectively,
\[ \tmr_p = \A \x - \tmb, \quad \tmr_d = \A^{\top} \y + \s - \tmc, \quad
   \tmr_{\mu} = \s - \mu \X^{- 1} \e . \]
Linearizing these relations, we obtain
\begin{align*}
  \A \Delta \x =  - \tmr_p, \quad
  \A^{\top} \Delta \y + \Delta \s = - \tmr_d,\quad
  \Delta \s + \mu \X^{- 2} \Delta \x =  - \tmr_{\mu},
\end{align*}
where $(\X + \Delta \X)^{- 1} \e$ is linearized as $\X^{- 1} \e -\X^{- 2} \Delta \x$. After eliminating $\Delta \x$ and $\Delta \s$, we find
\begin{align}
  \A \X^2 \A^{\top} \Delta \y ={} & -\mu \tmr_p + \A \X^2 ( \tmr_{\mu} -
  \tmr_d ) \nonumber\\
  \Delta \s ={} & - \tmr_d - \A^{\top} \Delta \y \label{eqn:infeaspkkt}\\
  \Delta \x ={} & - \tfrac{1}{\mu} \X^2 ( \tmr_{\mu} + \Delta \s ).
  \nonumber
\end{align}
Similar to primal-dual IPM, after obtaining $\Delta \y$ from the normal equation, we do backward
substitution to recover $\Delta \s$ and $\Delta \x$. Then, primal IPM performs a ratio test, decreases $\mu$, and goes to the next iteration. One advantage of \eqref{eqn:infeaspkkt} is that we explicitly maintain dual variables $(\y, \s)$, and primal IPM behaves like a primal-dual IPM. We conclude this section with a discussion of potential numerical issues arising from \eqref{eqn:infeaspkkt}: although mathematically correct, solving \eqref{eqn:infeaspkkt} directly
sometimes fails due to floating point error: error in $\Delta\x$ is magnified by $1/\mu$. 
%
To resolve this issue, we can rearrange and divide the normal equation by $\mu$ to obtain
\begin{align}
  \A \X^2 \A^{\top} ( \tfrac{1}{\mu} \Delta \y ) = & - \tmr_p +
  \tfrac{1}{\mu} \A \X ( \X \bs \e - \mu \e ) - \tfrac{1}{\mu} \A
  \X^2 \tmr_d, \nonumber
\end{align}

and solve for $\mu^{- 1} \Delta \y$ instead. Notice that $\tfrac{1}{\mu} ( \X \bs
\e - \mu \e )$ corresponds to the proximity measure and is bounded for primal IPM. Hence it is possible to obtain $\mu^{- 1} \Delta
\y$ with high accuracy so long as $\tmr_d$ is not large. This strategy provides a more
stable solution for the primal IPM direction.

\subsection{Dualization, Symmetric Standard Form and Variable Bounds}

Convergence of the normal matrix is only available for primal IPM, but our results extend to the dual problem in symmetric primal-dual standard form:
\begin{equation} \label{eqn:symstform}
\minf{\x \geq \0} \quad \langle \tmc, \x \rangle \quad 
\text{subject to} \quad \A \x \geq \tmb \quad \text{and} \quad
\maxf{\y \geq \0} \quad \langle \tmb, \y \rangle \quad 
\text{subject to} \quad \A^\top \y \leq \tmc.
\end{equation}
\eqref{eqn:symstform} ensures the primal and dual problems have the same form. 
Therefore, primal IPM can be applied after dualization. By introducing primal and dual slack $\w, \z \geq \0$ such that $\A \x - \w =
\tmb, \A^{\top} \y + \z = \tmc$,
\begin{align*}
	\A \x - \w  ={} \tmb, & \quad 
    ~\Z \X \e ={} \mu \e,  \quad 
    ( \x, \w )  \geq{} \0\\
    \A^{\top} \y + \z ={} \tmc, & \quad 
    \W \Y \e  ={} \mu \e, \quad
    ( \y, \z ) \geq{} \0
\end{align*}

In the primal form, $\z = \mu \X^{- 1} \e$ and $\y = \mu \W^{- 1}
\e$ are linearized:
\begin{align*}
  \Delta \z + \mu \X^{- 2} \Delta \x  ={}  \mu \X^{- 1} \e - \z \quad 
 \text{and} \quad  \Delta \y + \mu \W^{- 2} \Delta \w  ={}  \mu \W^{- 1} \e - \y .
\end{align*}
In this case, the normal matrix is $\A \X^2 \A^{\top} + \W^2$.
 In the dual form, $\x = \mu \Z^{- 1} \e, \w = \Y^{- 1} \e$ are linearized:
\begin{align*}
  \Delta \x + \mu \Z^{- 2} \Delta \z & ={}  \mu \Z^{- 1} \e - \x
  \quad  \text{and} \quad 
  \Delta \w + \mu \Y^{- 2} \Delta \y ={}  \mu \Y^{- 1} \e - \w .
\end{align*}
The normal matrix becomes $\A^{\top} \Y^2 \A + \Z^2$. The symmetric form allows switching between primal and dual problems to exploit sparsity or (non-)degeneracy. 

\paragraph{Variable Bound.} Another case of interest is the standard form with bounds,
\begin{equation}
  \minf{\x} \quad \langle \tmc, \x \rangle \quad
  \text{subject to} \quad \A \x = \tmb \quad
   \0 \leq \x \leq \tmu 
\end{equation}

where $u_j \in [0, \infty) \cup \{ + \infty \}$. The barrier function becomes
$- \log (u - x) - \log x$ for upper-bounded variables, and the barrier Hessian gives the scaling matrix $\D^2 = \X^{- 2} + ( \U - \X )^{-
2}$ with diagonals
\[ \sqrt{\tfrac{1}{x^{- 2} + (u - x)^{- 2}}} = \tfrac{x (u - x)}{\sqrt{x^2 +
   (u - x)^2}} . \]
The diagonal of the scaling matrix also remains stable when $x$ converges.

\subsection{Semidefinite Programming}

This section discusses how a preconditioned primal IPM generalizes to SDP. Consider standard form  SDP:
\[ \minf{\X} \quad \langle \mathbf{C}, \X \rangle \quad \text{subject
   to} \quad \mathcal{A} (\X) = \tmb, \quad \X \succeq \0, \]
where $\mathcal{A} (\X) = ( \langle \A_1, \X \rangle, \ldots,
\langle \A_m, \X \rangle )^\top \in \mathbb{R}^m$. The primal normal matrix is
\[ \M_{\X} = \left(\begin{array}{ccc}
     \left\langle \A_1, \X \A_1 \X \right\rangle & \cdots & \left\langle \A_1,
     \X \A_m \X \right\rangle\\
     \vdots & \ddots & \vdots\\
     \left\langle \A_m, \X \A_1 \X \right\rangle & \cdots & \left\langle \A_m,
     \X \A_m \X \right\rangle
   \end{array}\right) . \]
When $m$ is large, setting up and factorizing $\M_{\X}$ becomes a bottleneck in IPM
for SDP. Suppose $\M_{\X_k}$ is applied as a preconditioner. The action of the normal matrix $\M_{\X}$ on $\Delta \y$ can be written as 
\[ \M_{\X} \Delta \y =\mathcal{A} ( \X ( \mathcal{A}^{\ast} \Delta
   \y ) \X ). \]

Hence, we can solve the normal equation at each iterate $\x$ using an iterative solver without constructing $\M_{\X}$ explicitly. Note that forming the preconditioner still requires (occasional) explicit construction and factorization of the primal normal matrix.

\subsection{Combining the Benefits of Primal IPM with Primal-Dual IPM}
We conclude with a few remarks on how primal IPM can benefit existing primal-dual solvers.  First of all, we advocate that primal IPM should not be used as a replacement for primal-dual IPM but rather as a complement: 
the advantage of primal IPM only appears near convergence, while the initial iterations needed for primal IPM to ensure proximity in the (thresholded) distance can be slower than the time required for primal-dual to solve the whole problem. 
Therefore, it works best to let primal-dual IPM iterate until a medium-accuracy solution is found and then to invoke primal IPM to complete the final iterations. 
Last, we discuss the complexity of implementation. 
On an algorithmic level, 
infeasible-start primal IPM is similar to primal-dual IPM; the main difference is simply the use of an iterative solver. 
However, given that primal-dual IPM software often uses an iterative refinement procedure (a particular type of iterative method),  
primal IPM can be implemented in a primal-dual solver without much additional effort.

\section{Numerical Experiments} \label{sec:exp}

In this section, we conduct numerical experiments to validate the performance of primal IPM.

\subsection{Experiment Setup}
Our experiment contains three parts.
\begin{itemize}[leftmargin=10pt,itemsep=2pt]
\item \tmtextbf{Section \ref{exp:lp}} demonstrates the practical utility of a combined primal-dual and primal IPM  for solving LPs.

  \item \tmtextbf{Section \ref{exp:sdp}} demonstrates the effectiveness of primal IPM for SDP with many constraints.
  
  \item \tmtextbf{Section \ref{exp:conv}} studies the convergence behavior of
  $\{ \x_k \}$ in primal-dual IPM.
\end{itemize}

\paragraph{Testing environment.} Qualitative results in \tmtextbf{Section \ref{exp:conv}} are obtained with \texttt{Macbook Pro} with
\texttt{Apple Silicon 32G}. The rest of the
results are obtained with \texttt{AMD Ryzen 9 5900X 12-Core 128G}.

\paragraph{Dataset.} In LP-related tests, we use 924 LP instances from the LP
relaxation of MIPLIB 2017 \cite{miplib2017}. Each instance is
presolved by \texttt{Gurobi} and converted to standard form. We use the instances from Hans Mittelmann's SDP benchmark \cite{mittelmann2020benchmarking} for SDP-related tests.

\paragraph{IPM solver.} We compare three primal-dual IPM implementations:
\begin{itemize}[leftmargin=10pt, itemsep=2pt]
  \item \texttt{Gurobi} and \texttt{COPT} for IPM convergence
  behavior analysis in \tmtextbf{Section \ref{exp:conv}}.
  
  \item A textbook version of a primal-dual IPM with Mehrotra type
  corrector step in \texttt{Matlab} for {\tmstrong{Section \ref{exp:conv}}}.
  
  \item A textbook version of a primal-dual IPM within the open-source
  software \texttt{HDSDP} \cite{gao2022hdsdp} for \tmtextbf{Section \ref{exp:lp}}.
\end{itemize}
We compare two primal IPM implementations:
\begin{itemize}[leftmargin=10pt, itemsep=2pt]
  \item The same primal IPM as in \texttt{HDSDP}. We use either
  conjugate gradient method as the iterative solver. This
  version is used for \tmtextbf{Section \ref{exp:lp}}.
  
  \item SDP version of primal IPM in \texttt{HDSDP}. The conjugate gradient method is used as the iterative solver. This version is used for
  \tmtextbf{Section \ref{exp:sdp}}.
\end{itemize}
\paragraph{Convergence metric.} For LPs, our experiments use standard convergence metrics: primal infeasibility, dual infeasibility, and duality gap:
\[ e_p \assign \tfrac{\| \A \x - \tmb \|}{1 + \| \tmb
   \|}, \quad e_d \assign \tfrac{\| \A^{\top} \y + \s - \tmc
   \|}{1 + \| \tmc \|}, \quad e_g \assign \tfrac{|
   \langle \tmc, \x \rangle - \langle \tmb, \y \rangle
   |}{1 + | \langle \tmc, \x \rangle | + |
   \langle \tmb, \y \rangle |}. \]
For SDPs, our experiments use DIMACS error \cite{mittelmann2003independent}. 
\subsection{Primal IPM for LPs} \label{exp:lp}

In this section, we report the performance of primal IPM applied to solving LPs.

\paragraph{Additional experiment details.} For each LP instance, we run 
\begin{itemize}[leftmargin=10pt]
	\item  Standard primal-dual IPM with predictor-corrector step.
	\item  Hybrid primal-dual-primal IPM. The method starts with primal-dual IPM. When primal-dual IPM runs, the solver estimates the ratio between the time factorizing the normal matrix and forward/backward substitution. The hybrid solver switches to primal IPM if $\|\x_{k} - \x_{k-1}\|_{\x_k, \nu} \leq 10^{-1}$ and the time ratio is greater than 30.
For convenience, we call this hybrid algorithm primal IPM.
\end{itemize}
until one of the three conditions is satisfied: \textbf{1)} $\max \{e_p, e_d, e_g \} \leq 10^{-10}$, in which case we say the solver succeeds on the instance; \textbf{2)} the iteration number reaches the maximum limit of 100; \textbf{3)} matrix factorization fails due to insufficient memory or numerical error. In the latter two cases, we say the solver fails the test on the instance. There are four possible cases:

\begin{enumerate}[leftmargin=25pt, label={\textbf{C\arabic{*}}.},ref={\textbf{C\arabic{*}}}]
	\item Primal-dual IPM fails the test.\label{case:1}
	\item The hybrid solver does not invoke Primal IPM. \label{case:2}
	\item Primal-dual IPM passes the test, while primal IPM fails the test. \label{case:3}
	\item Both primal-dual IPM and primal IPM pass the test. \label{case:4}
\end{enumerate}
Let $T_{PD}$ and $T_{P}$ be the solution time of primal-dual and primal IPM in some instance, respectively. We define the speedup of primal IPM by $s = \frac{T_{PD} - T_{P}}{T_{PD}}$. If $T_{PD} \leq 15$, we repeat the experiment 10 times and record the average running time. Detailed time statistics are given in the \textbf{Appendix \ref{app:stats}}.

\begin{figure}[h]
\centering
\includegraphics[width=0.4\textwidth]{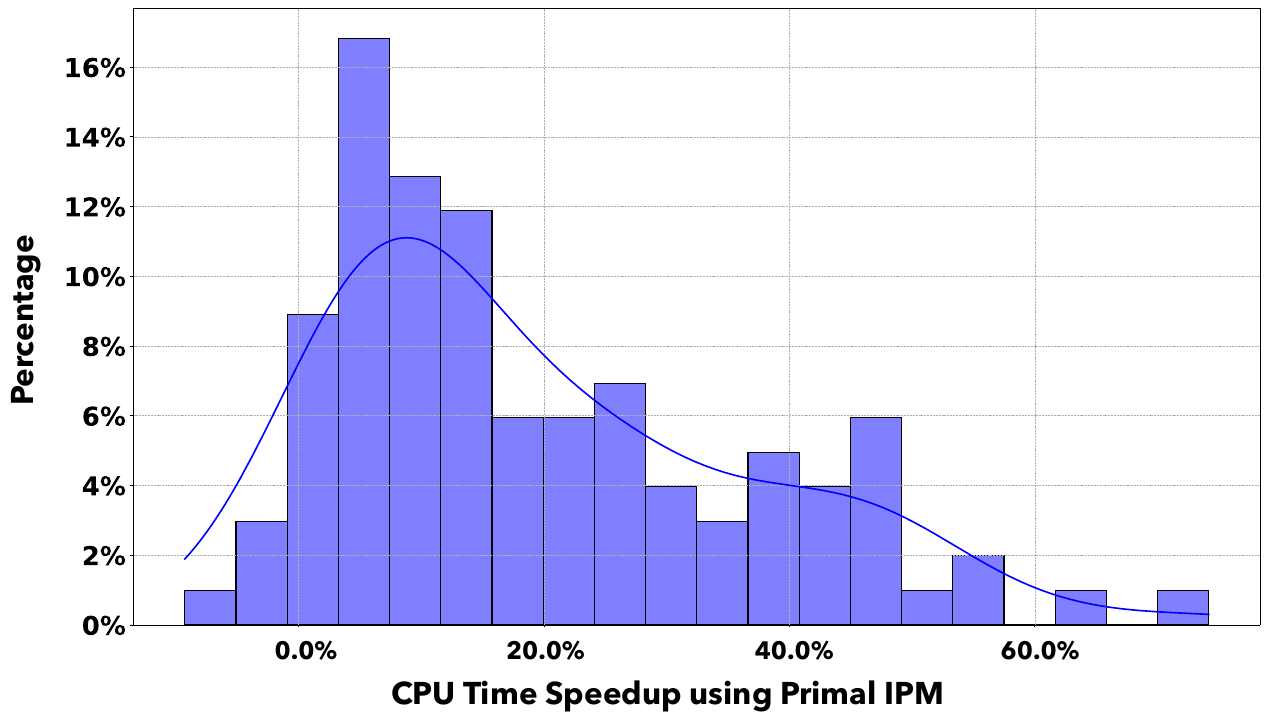}
\includegraphics[scale=0.25]{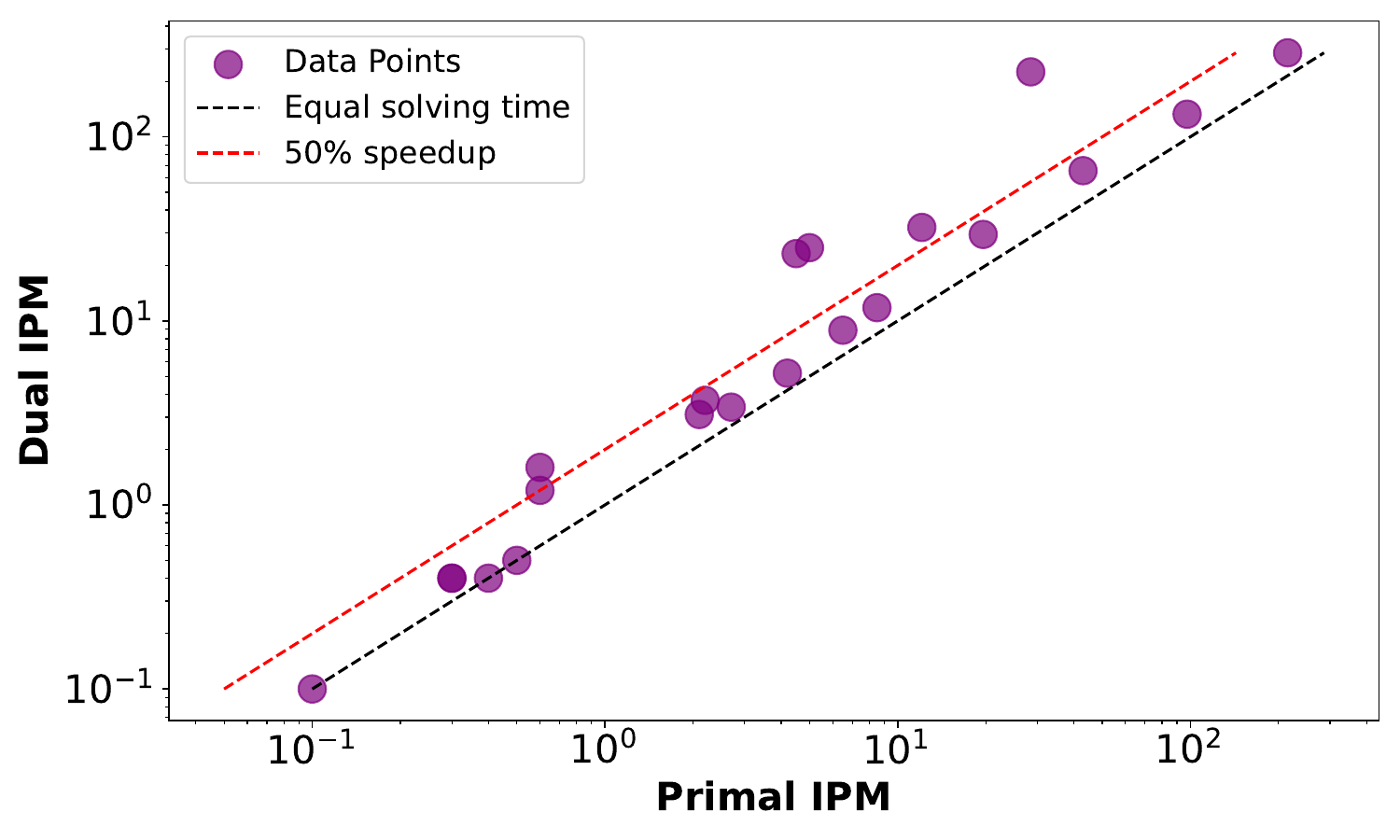}
\caption{Left: distribution of primal IPM speedup over 139 LPs. Right: Dual vs hybrid dual-primal IPM on 21 SDPs. The hybrid approach is always faster, with at least 50\% speedup on 6 instances.
\label{fig:pipm-distri}}
\end{figure}

\paragraph{Experiment results. } With the setup above, 744 instances fall into \ref{case:1} and \ref{case:2}. Among the remaining 180 instances, 41 are in \ref{case:3} and 139 are in \ref{case:4}. Primal IPM shows benefits in \ref{case:4}. \textbf{Figure \ref{fig:pipm-distri}} (Left) plots the distribution of speedup from primal IPM in \ref{case:4}. Only six instances slow down when primal IPM is invoked and passes the test, and more than $30\%$ of instances speed up by at least $30\%$. Around $10\%$ of instances speed up by at least $50\%$. This result suggests primal IPM can be useful for a broad class of LPs.

\subsection{Primal IPM for SDPs} \label{exp:sdp}
We implement preconditioned primal IPM in \texttt{HDSDP} (called \texttt{PSDP}) and test on 21 instances with $m \geq 10n$ from the SDP benchmark \cite{mittelmann2020benchmarking} and the DIMACS challenge \cite{mittelmann2003independent}. Both \texttt{HDSDP} and \texttt{PSDP} solve to $10^{-6}$  accuracy in DIMACS error. \textbf{Figure \ref{fig:pipm-distri}} (Right) illustrates the speed up using primal IPM when it is combined with a dual IPM method.
When $m \gg n$, primal IPM significantly reduces the solution time for the SDP solver (up to 10 times faster). This observation aligns with our intuition since, in this case, most of the solution time is spent on forming and factorizing the normal matrix. Moreover, since SDP generally does not have superlinear convergence, the potential speedup due to switching to the primal IPM is larger.

\subsection{Convergence Behavior of the Primal-Dual IPM Iteration Sequence} \label{exp:conv}
\begin{figure}[h]
\centering
\includegraphics[scale=0.17]{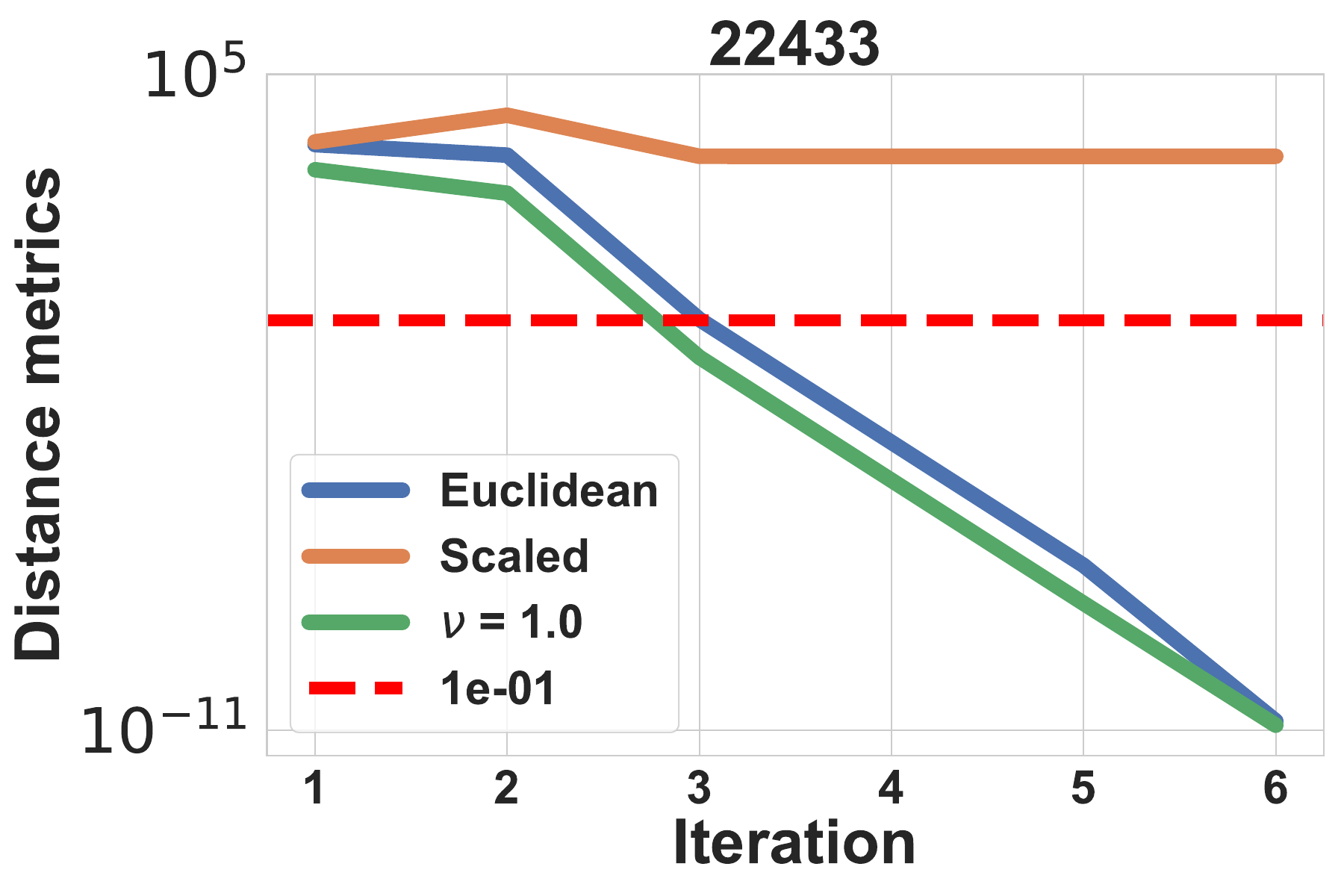}
	\includegraphics[scale=0.17]{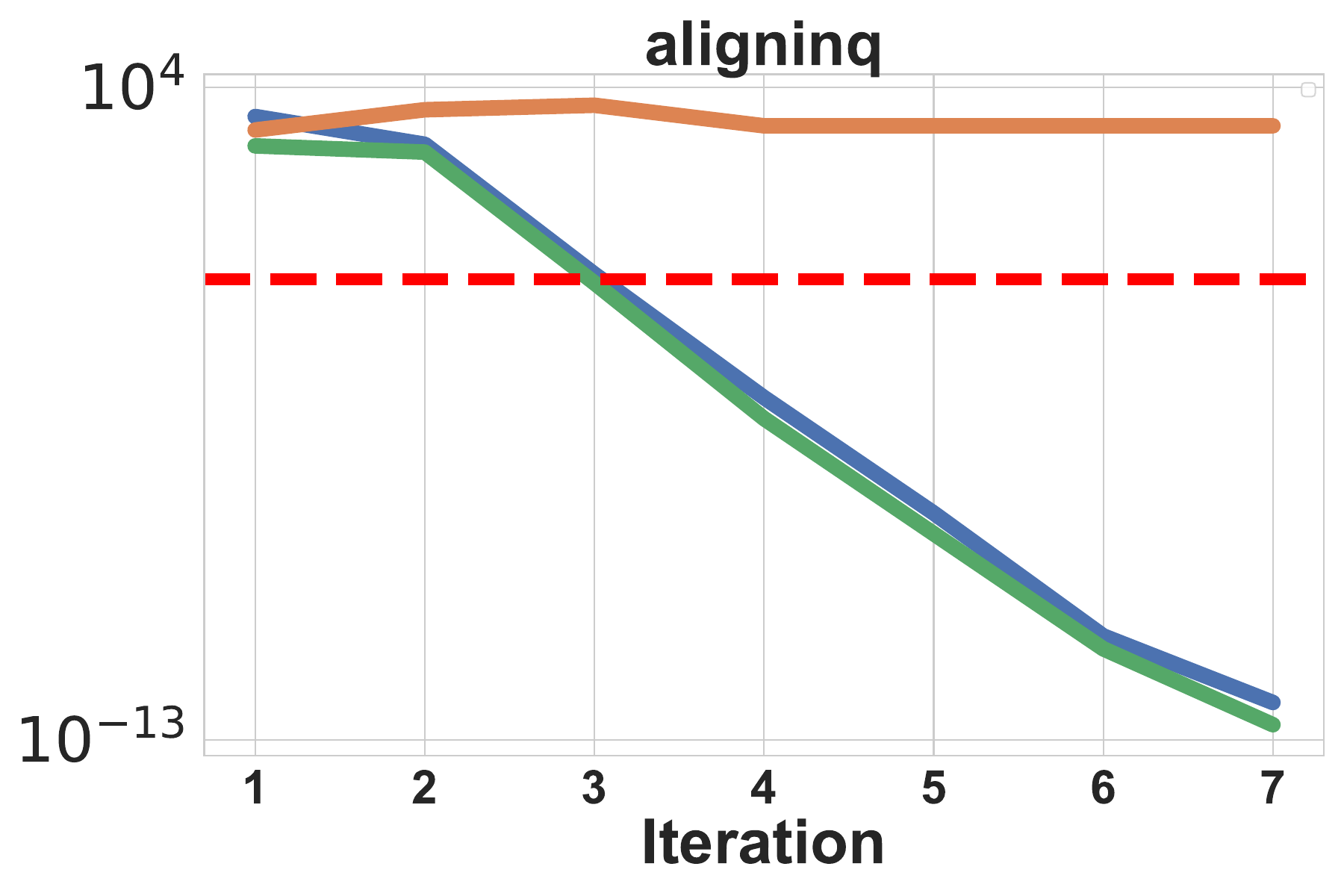}\includegraphics[scale=0.17]{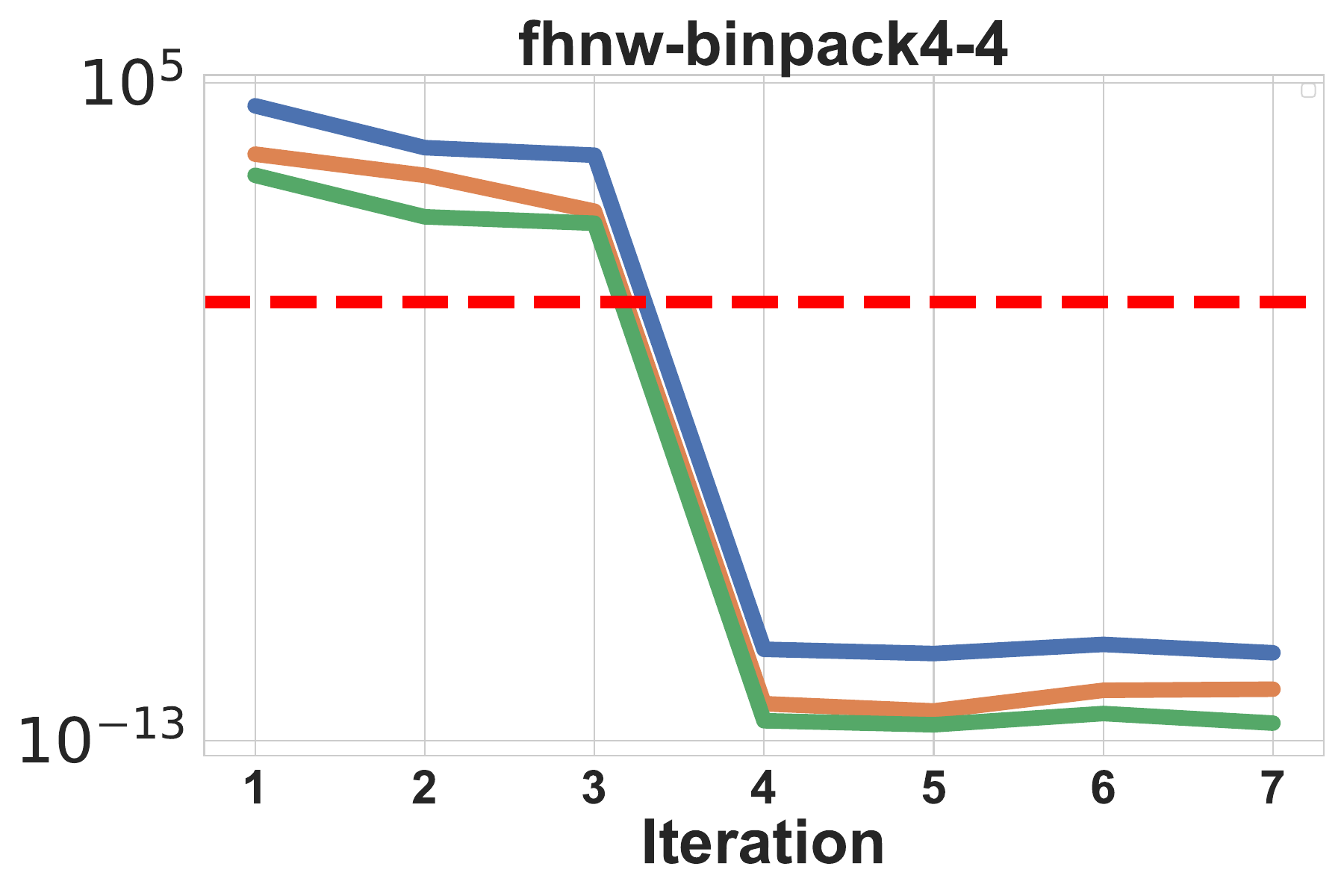}
	\caption{Problems with convergence pattern: linear convergence \label{fig:linconv}}
\end{figure}

The convergence behavior of primal-dual IPM determines the performance improvement from switching to primal IPM: switching to primal IPM offers the largest benefit when primal-dual IPM converges slowly in the last iterations. Therefore, in this section, we discuss the convergence behavior of the primal-dual IPM iteration sequence. We take 924 instances from the {MIPLIB} and record $\{ \|
\x_{k + 1} - \x_k \| \}, \{ \| \x_{k + 1} - \x_k
\|_{\x_k, 1} \}$ from primal-dual.  We summarize the general convergence patterns of IPM as follows.
\paragraph{Consistent linear convergence (with plateau).} All distance metrics except scaled distance decrease linearly (\textbf{Figure \ref{fig:linconv}}) since early stage of IPM. There is no noticeable superlinear convergence behavior. IPM iterations may or may not reach a plateau region. Problems with zero objective $\tmc = \0$ are more likely to exhibit this type of convergence pattern. Primal IPM speeds up convergence by 30\% to 60\% on these problems.

\begin{figure}[h]
\centering
	\includegraphics[scale=0.17]{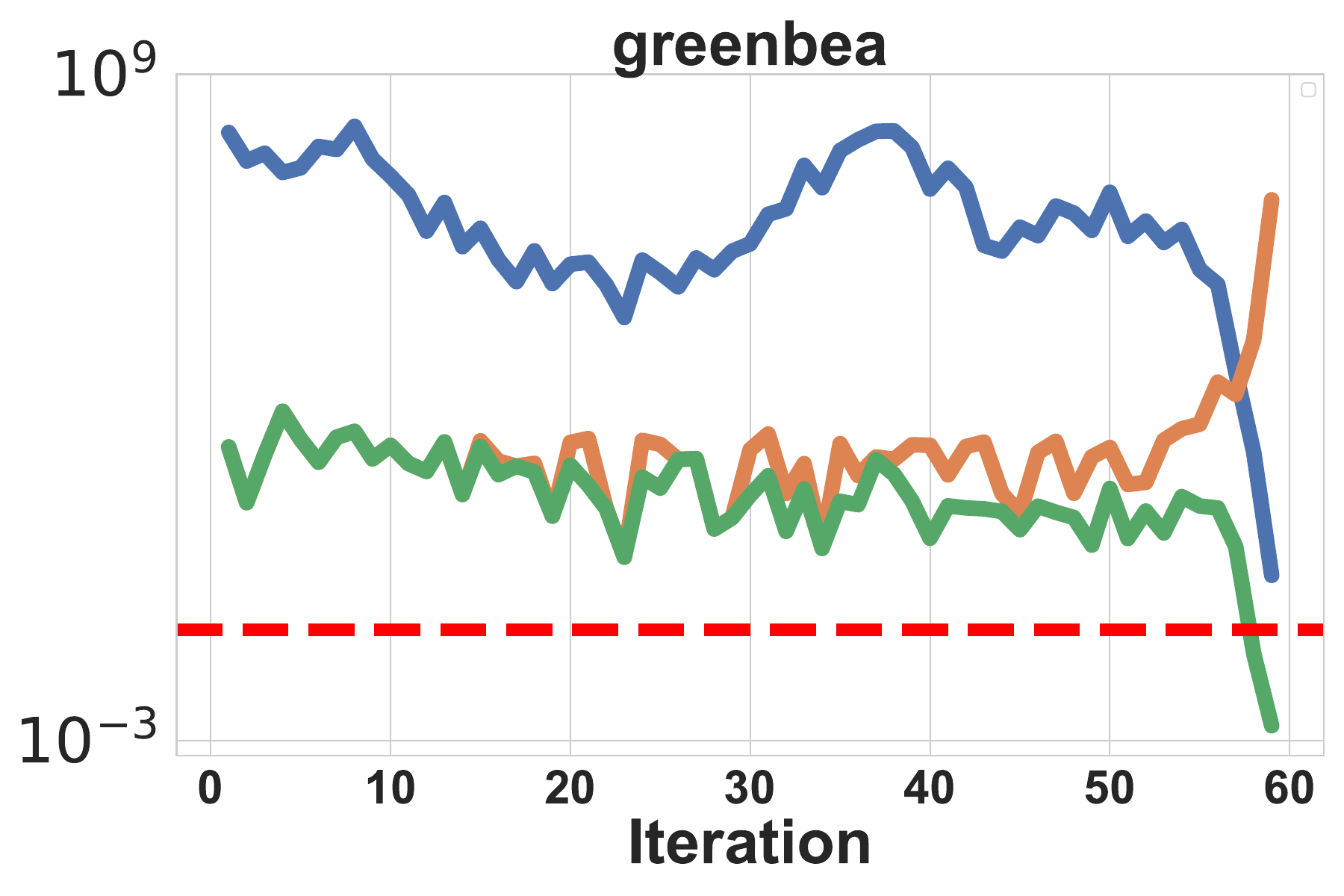}\includegraphics[scale=0.17]{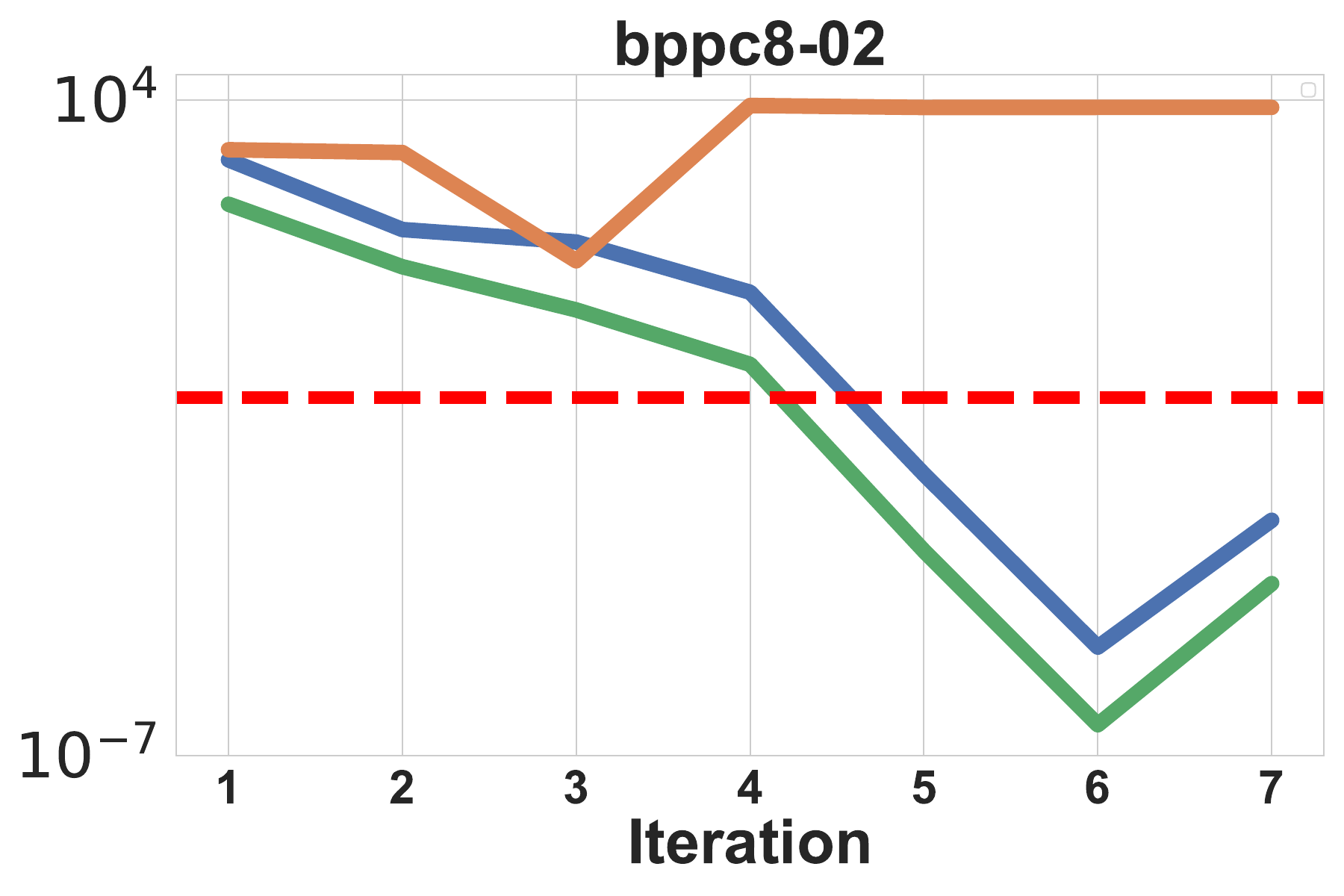}
	\includegraphics[scale=0.17]{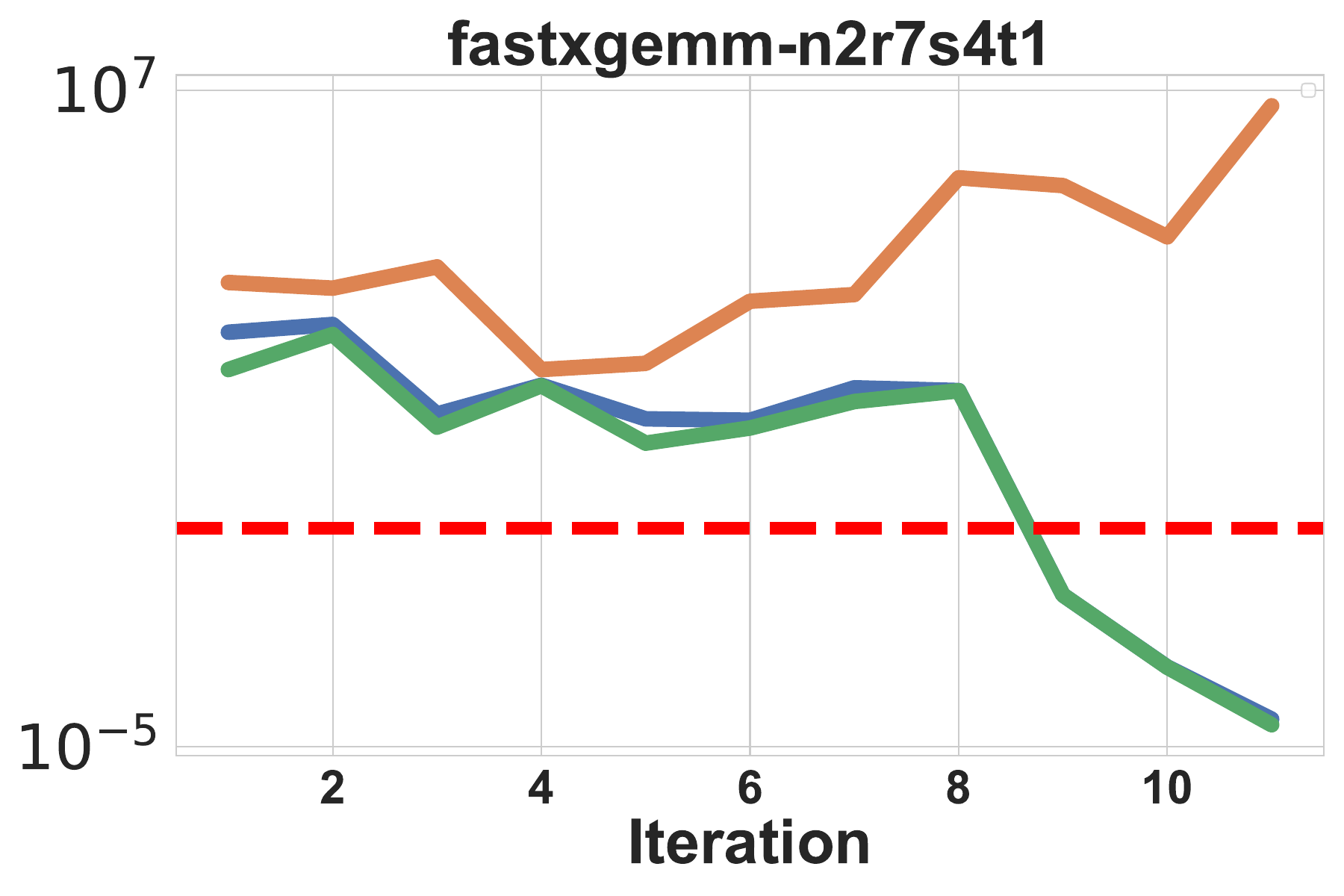}
	\caption{Problems with convergence pattern: fast convergence at the end\label{fig:fastend}. For problems like the leftmost figure, primal IPM  is not useful. In contrast, for the second and third figures, the number of iterations after entering the fast convergence phase is comparable to the total iteration count, so primal IPM offers a speedup.}
\end{figure}

\paragraph{Rapid convergence only at the end.}\textbf{Figure \ref{fig:fastend}} shows that for some problems, all the distance metrics decrease at the slow linear or sublinear rate in the initial iterations of IPM. In the last 2-5 iterations, (super)linear convergence appears, and iterates converge fast. More than half of the tested problems exhibit this convergence pattern. Whether primal IPM hastens convergence is determined by \textbf{1)} the convergence speed of thresholded distance; \textbf{2)} the duration of the superlinear convergence phase; \textbf{3)} the cost of factorizing the normal equations compared to forward/backward solve. The speedup from primal IPM is between 5\% to 15\% in this case. There is no benefit in applying primal IPM if superlinear convergence reaches the optimal solution in one or two iterations. 

\begin{figure}[h]
\centering	\includegraphics[scale=0.17]{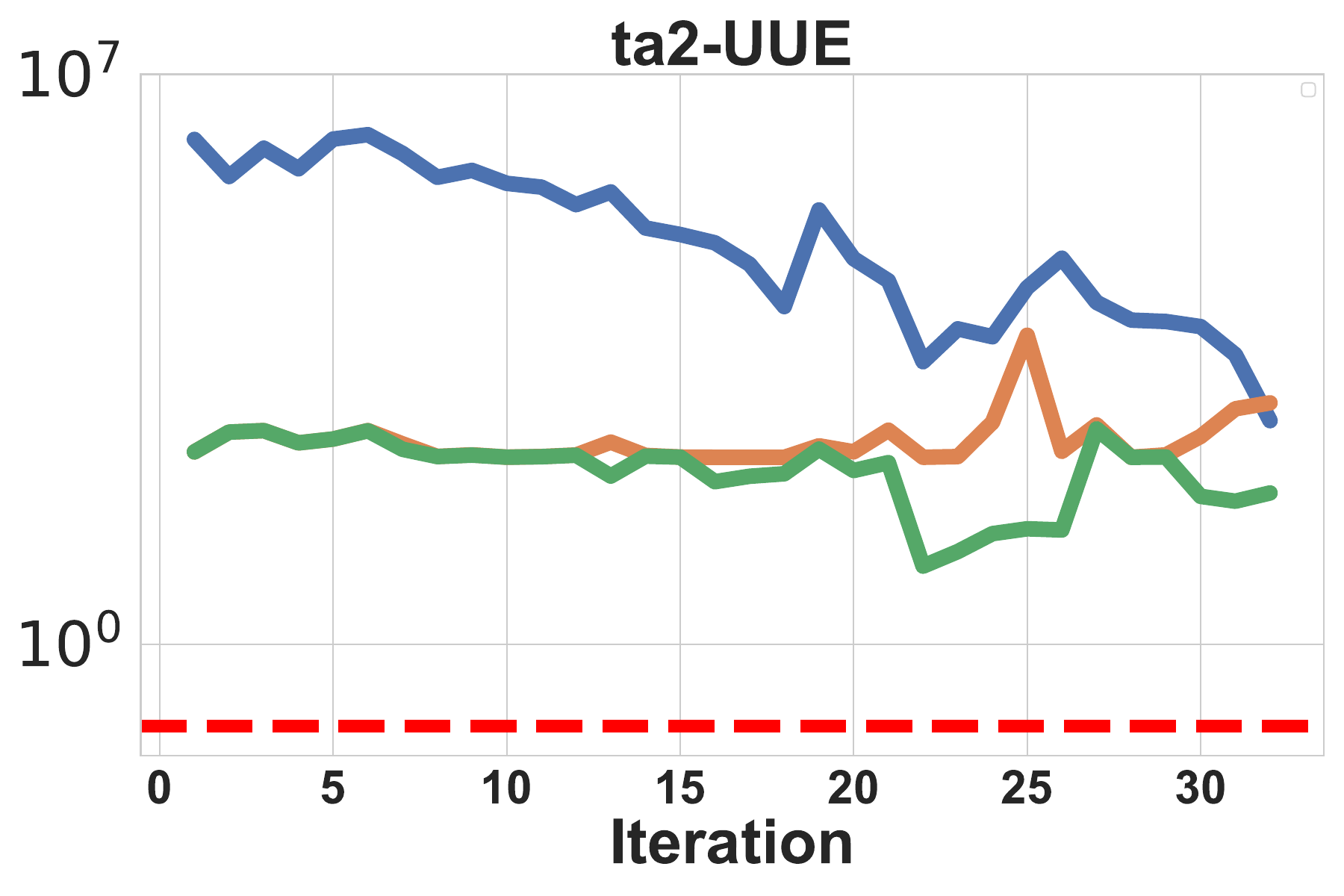}\includegraphics[scale=0.17]{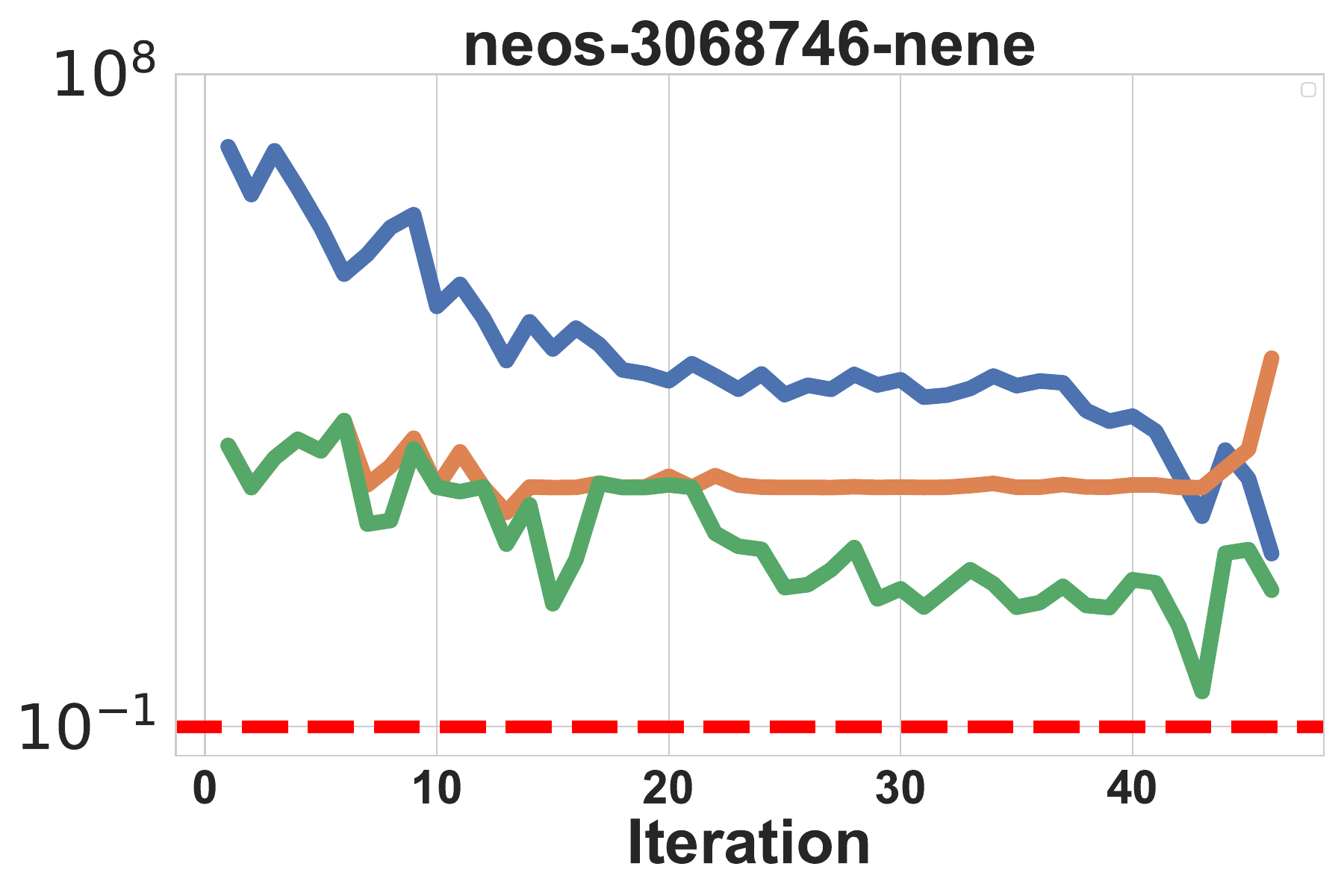}		\includegraphics[scale=0.17]{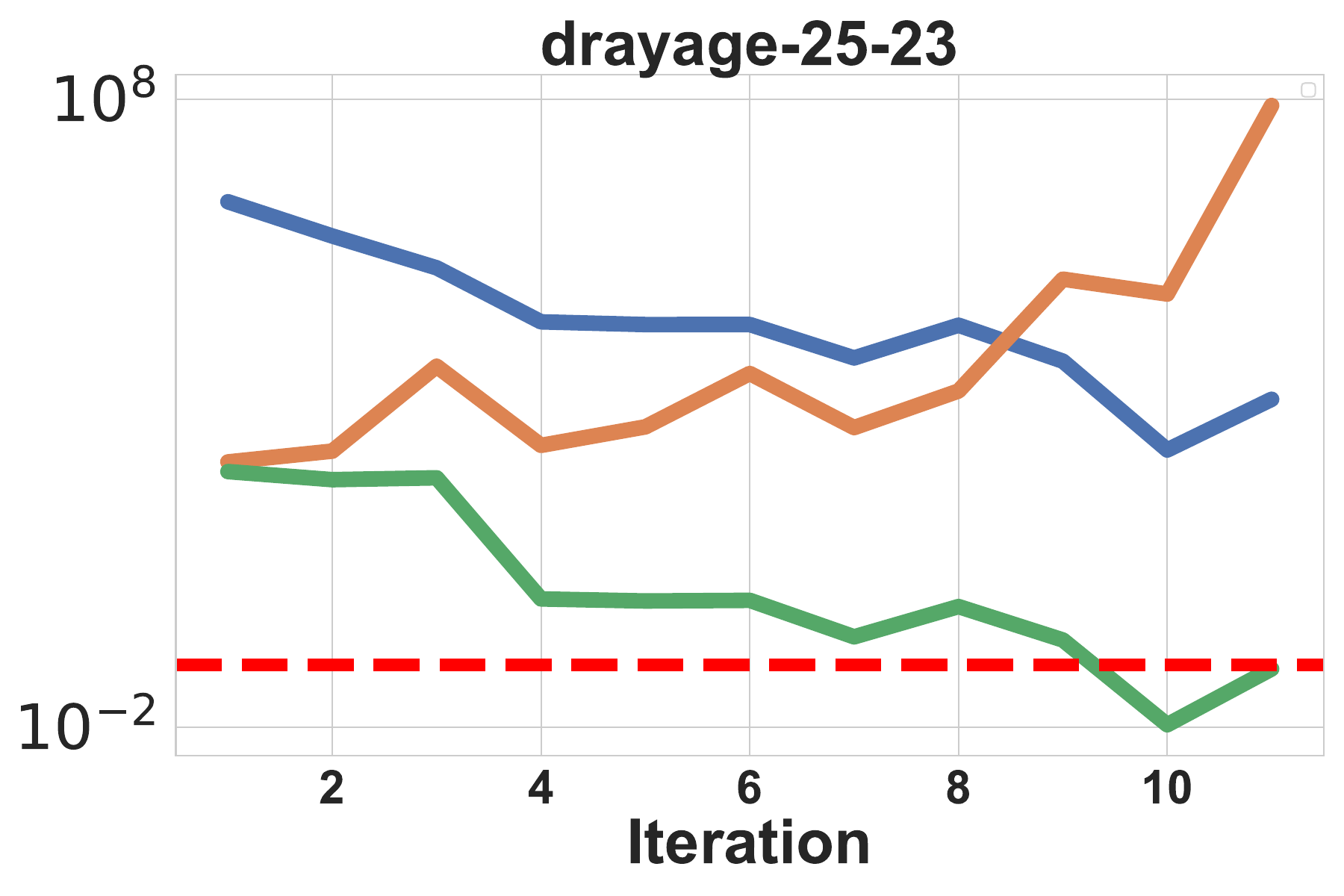}
	\caption{Problems with convergence pattern: no clear convergence. \label{fig:noclearconv} Often this pattern results when the solution $\x$ has extremely imbalanced coordinates. }
\end{figure}

\paragraph{No clear convergence.} For the rest of the problems (less than 5\%), like the examples shown in \textbf{Figure \ref{fig:noclearconv}}, the iteration sequence does not indicate clear convergence behavior. This case is less common and generally often happens when the primal solution has large-scale elements: the objective value converges while the solution sequence oscillates. Our theory excludes the use of primal IPM for such problems.

\section{Conclusions}\label{sec:conclusion}

In this paper, we revisit the primal interior point method and explore its potential advantages over the widely-used primal-dual IPM. We leverage the stability of the primal normal matrix to enhance computational efficiency in the primal IPM, particularly as the algorithm approaches convergence. Extensive numerical experiments demonstrate the practical efficiency of primal IPM on a broad class of LPs.

\section*{Acknowledgement}

The authors are grateful to Dr. Joachim Dahl and Dr. Christian Bliek from the COPT development team, and Zachary Frangella for their discussion and feedback on this work. MU gratefully acknowledges support from the National Science Foundation (NSF) Award IIS-2233762, the Office of Naval Research (ONR) Award N000142212825 and N000142312203, and the Alfred P. Sloan Foundation.

\renewcommand \thepart{}
\renewcommand \partname{}

\bibliography{ref.bib}
\bibliographystyle{plain}

\doparttoc
\faketableofcontents
\part{}

\newpage
\appendix

\addcontentsline{toc}{section}{Appendix}
\part{Appendix} 
\parttoc

\paragraph{Structure of the Appendix} 

The appendix is organized as follows. In \textbf{Section \ref{app:proof-ipmscal}}, we prove the main results in our paper. \textbf{Section \ref{app:inexact}} is a self-contained part that leverages our main results to establish convergence of inexact primal IPM.
\newpage
\begin{landscape}
\section{Detailed statistics on MIPLIB instances} \label{app:stats}
\begin{table}[h]
\centering
\resizebox{1.096\textwidth}{!}{
  \begin{tabular}{cccccccc}
    \toprule
    Instance & Primal Time & Primal-dual Time & Speedup & Instance & Primal
    Time & Primal-dual Time & Speedup\\
    \midrule
    allcolor58 & 1003.678 & 1045.208 & 0.041378 & neos-4338804-snowy & 0.187 &
    0.223 & 0.193895\\
    berlin\_5\_8\_0 & 0.182 & 0.224 & 0.231925 & neos-4797081-pakoka & 1.585 &
    2.116 & 0.334746\\
    bnatt400 & 0.505 & 0.619 & 0.225389 & neos-4971100-bodrog & 5.441 & 5.660
    & 0.040281\\
    bnatt500 & 0.911 & 1.120 & 0.229308 & neos-498623 & 0.574 & 0.587 &
    0.023318\\
    cdc7-4-3-2 & 144.080 & 164.713 & 0.143203 & neos-873061 & 148.944 &
    153.862 & 0.033016\\
    chromaticindex128-5 & 13.771 & 15.711 & 0.140903 & neos-885524 & 1.579 &
    1.433 & -0.092729\\
    chromaticindex256-8 & 103.370 & 118.173 & 0.143205 & neos-950242 &
    2034.914 & 2143.839 & 0.053528\\
    chromaticindex32-8 & 0.220 & 0.280 & 0.273078 & neos-960392 & 7.323 &
    9.232 & 0.260649\\
    chromaticindex512-7 & 805.685 & 922.078 & 0.144465 & neos9 & 64.425 &
    74.273 & 0.152860\\
    CMS750\_4 & 19.897 & 24.745 & 0.243658 & nexp-150-20-8-5 & 1.208 & 1.254 &
    0.037911\\
    cod105 & 0.113 & 0.160 & 0.413499 & no-ip-64999 & 0.273 & 0.386 &
    0.415453\\
    cryptanalysiskb128n5obj14 & 27.533 & 38.303 & 0.391161 & no-ip-65059 &
    0.260 & 0.371 & 0.426903\\
    cryptanalysiskb128n5obj16 & 28.689 & 39.961 & 0.392899 & ns1208400 & 0.361
    & 0.449 & 0.244861\\
    cvrpa-n64k9vrpi & 431.431 & 479.280 & 0.110908 & ns1828997 & 3332.151 &
    4068.890 & 0.221100\\
    cvrpb-n45k5vrpi & 51.542 & 56.801 & 0.102025 & osorio-cta & 0.608 & 0.685
    & 0.125941\\
    dale-cta & 19.820 & 20.415 & 0.030033 & probportfolio & 0.028 & 0.030 &
    0.067266\\
    disctom & 0.073 & 0.101 & 0.377128 & pythago7824 & 4.204 & 6.261 &
    0.489257\\
    genus-g31-8 & 492.335 & 601.367 & 0.221459 & pythago7825 & 4.618 & 6.877 &
    0.489196\\
    genus-g61-25 & 1050.716 & 1311.999 & 0.248671 & queens-30 & 0.197 & 0.224
    & 0.139101\\
    graph20-80-1rand & 7.089 & 7.496 & 0.057381 & ramos3 & 0.494 & 0.656 &
    0.327558\\
    graph40-20-1rand & 68.218 & 73.883 & 0.083047 & scpj4scip & 1.025 & 1.087
    & 0.060606\\
    graph40-40-1rand & 958.515 & 956.266 & -0.002346 & scpk4 & 2.280 & 2.262 &
    -0.008056\\
    lectsched-2 & 2.234 & 2.859 & 0.279880 & shipsched & 2.680 & 3.053 &
    0.139111\\
    lectsched-4-obj & 0.572 & 0.606 & 0.059112 & shipschedule6shipsmixi &
    0.446 & 0.599 & 0.343217\\
    lectsched-5-obj & 11.697 & 13.078 & 0.118042 & shipschedule8shipsmixuci &
    0.547 & 0.753 & 0.376145\\
    mspsphard01i & 0.356 & 0.467 & 0.310907 & supportcase18 & 0.191 & 0.228 &
    0.195104\\
    neos-1324574 & 0.555 & 0.606 & 0.092235 & supportcase21i & 60.069 & 66.351
    & 0.104575\\
    neos-2978193-inde & 1.100 & 1.204 & 0.094842 & supportcase29 & 7.555 &
    7.529 & -0.003391\\
    neos-2978205-isar & 97.132 & 97.086 & -0.000471 & supportcase30 & 0.057 &
    0.084 & 0.471918\\
    neos-3004026-krka & 13.471 & 20.117 & 0.493343 & supportcase3 & 2.014 &
    2.991 & 0.485151\\
    neos-3135526-osun & 0.655 & 0.715 & 0.091298 & supportcase4 & 0.334 &
    0.495 & 0.481269\\
    neos-3211096-shag & 1.335 & 1.903 & 0.425468 & traininstance2 & 1.520 &
    1.594 & 0.048841\\
    neos-3355323-arnon & 13.661 & 21.152 & 0.548381 & tw-myciel4 & 5.853 &
    6.203 & 0.059735\\
    neos-3426132-dieze & 0.171 & 0.203 & 0.189664 & v150d30-2hopcds & 46.815 &
    50.653 & 0.081980\\
    neos-3572885-glina & 7.999 & 9.097 & 0.137252 & van & 20.796 & 22.234 &
    0.069160\\
    neos-3603137-hoteo & 2.329 & 3.593 & 0.542789 &  &  &  & \\
    \bottomrule
  \end{tabular}
  }
  \caption{Speedup on instances where primal IPM is triggered}
\end{table}
	
\end{landscape}
\newpage
\section{Proof of Results in Section \ref{sec:pipm}} \label{app:proof-ipmscal}

\subsection{Auxiliary Results}
This section uses the following auxiliary results on primal central path solution $\x^\mu$.

\begin{lem}[\cite{roos1997theory}] \label{lem:pathproxbnd}
Under \ref{A1}, if $\delta ( \x, \mu ) \leq \lambda \leq
  \frac{\sqrt{2}}{2}$, then 
\begin{align}\label{ineq:dist-to-central}
    \| \x -
  \x^{\mu} \|_{\x^\mu} \leq \sqrt{1 - \sqrt{1 - 2 \lambda^2}}.
\end{align}
\end{lem}

\begin{lem} \label{lem:cplb}
  Under \ref{A1}, there exist $\alpha, \gamma > 0$ such that
  \begin{align}
    \gamma \mu \geq{}  x_j^{\mu},  \text{~for all}~~  j \in \mathcal{N}, \nonumber
  \end{align}
where $\Bcal, \Ncal$ is the maximal complementarity pair from \textbf{Definition \ref{def:maxcomp}}.
\end{lem}

\begin{proof}
Let $(\x^{\star}$, $\s^{\star})$ be a strictly complementary solution. By definition,
  \[ \x_{\mathcal{B}}^{\star} > \0 \quad \text{and} \quad  \s_{\mathcal{N}}^{\star} > \0 . \]
Given any $\x^{\mu} \in \mathcal{C}_p$, , there exists some dual feasible solution $(\s^\mu, \y^\mu)$ such that $ \X^{\mu} \s^{\mu} = \mu \e $ \cite{roos1992polynomial}. Then,
     \[
     \langle \x^{\mu} - \x^{\star}, \s^{\mu} - \s^{\star} \rangle = \langle \x^{\mu} - \x^{\star}, \A^\top (\y^{\star} - \y^{\mu}) 
  \rangle = \langle \A (\x^{\mu} - \x^{\star}), \y^{\star} - \y^{\mu} \rangle = 0,
     \]
where the first equality holds since $\s^{\mu} = \tmc - \A^{\top} \y^{\mu}$ and $\s^{\star} = \tmc -
     \A^{\top} \y^{\star}$, and the last equality holds since $\A \x^\star = \A \x^\mu = \tmb$.
  On the other hand,
  \begin{align}
    0 = \langle \x^{\mu} - \x^{\star}, \s^{\mu} - \s^{\star} \rangle =  \mu n - \langle \x^{\mu}, \s^{\star} \rangle - \langle
    \s^{\mu}, \x^{\star} \rangle
    = \mu n - \langle \x^{\mu}_{\mathcal{N}}, \s^{\star}_{\mathcal{N}}
    \rangle - \langle \s^{\mu}_{\mathcal{B}},
    \x^{\star}_{\mathcal{B}} \rangle \label{eqn:proof-2.2-2},
  \end{align}
  where the second equality holds due to the fact that $ \X^{\mu} \s^{\mu} = \mu \e $ and $\langle \x^\star, \s^\star \rangle = 0$, and the last one holds since  $\x^\star_\Ncal = \0$ and  $\s^\star_\Bcal = \0$.
Next we rearrange the terms to get  $\langle \x^{\mu}_{\mathcal{N}}, \s^{\star}_{\mathcal{N}}
     \rangle + \langle \s^{\mu}_{\mathcal{B}},
     \x^{\star}_{\mathcal{B}} \rangle = \mu n.$
  The nonnegativity of all the entries implies that $x_j^{\mu} s^{\star}_j \leq \mu n$ for all $j \in \mathcal{N}$, which implies
  \begin{align*}
    x_j^{\mu} \leq{}  \mu n (s_j^{\star})^{- 1} \leq \mu n ( \|
    \s^{\star}_{\mathcal{N}} \|_{- \infty} )^{- 1}, \, \text{~for all~} j \in
    \mathcal{N}.
  \end{align*}
Using \eqref{eqn:gamma}, $\gamma =  \frac{n}{\| \s^{\star}_{\mathcal{N}} \|_{-
  \infty}}$,
this completes the proof.
\end{proof}

\begin{lem} \label{lem:app:cplb}
  Given $\x \in \mathcal{F}_p^0, \mu \in (0, 1]$ such that $\delta ( \x,
  \mu ) \leq \lambda = \frac{1}{2}$, there exists some $\alpha > 0$ such
  that
  \[ \|\x\|_{-\infty} \geq \alpha \mu \]
\end{lem}

\begin{proof}
Again we use
  $\x^{\mu}, \s^{\mu}$ to denote primal and dual solutions on the central path with barrier parameter $\mu$.
  Now take any $\x' \in \mathcal{F}_p^0, \s' \in \mathcal{F}_d^0$ and
  \[ 0 = \langle \x^{\mu} - \x', \s^{\mu} - \s' \rangle = \mu n -
     \langle \x^{\mu}, \s' \rangle - \langle \s^{\mu}, \x'
     \rangle + \langle \x', \s' \rangle . \]
  Next rearrangement gives $\langle \s^{\mu}, \x' \rangle +
  \langle \x^{\mu}, \s' \rangle \leq \mu n + \langle \x', \s'
  \rangle$ and
  \[ \mu (x_j^{\mu})^{- 1} x'_j = s^{\mu}_j x'_j \leq \langle \s^{\mu},
     \x' \rangle + \langle \x^{\mu}, \s' \rangle \leq \mu n +
     \langle \x', \s' \rangle, \]
for all $j$,  which implies $x_j^{\mu} \geq \frac{\mu x'_j}{ \mu n + \langle
  \x', \s' \rangle } \geq \frac{\mu x'_j}{n + \langle \x',
  \s' \rangle}$ since $\mu \leq 1$. Now take some $\x \in
  \mathcal{F}_p^0$ and by \textbf{Lemma \ref{lem:pathproxbnd}},
  \begin{align}
    \max_{1 \leq j \leq n}  \Big| \tfrac{x_j - x^{\mu}_j}{x^{\mu}_j} \Big|
    ={} & \| ( \X^{\mu} )^{- 1} ( \x - \x^{\mu} )
    \|_{\infty} 
    \leq{}  \| \x - \x^{\mu}
    \|_{\x^\mu} \leq{} \sqrt{1 - \sqrt{1-\tfrac{1}{2}}}
    < 0.9\nonumber.
  \end{align}
  
  Hence for each $j$,
  \[ x_j \geq (1-0.9)\cdot x^{\mu}_j
     \geq \tfrac{0.1x'_j}{n + \langle \x', \s' \rangle}  \mu \]
  and taking $\alpha = \frac{0.1\|\x'\|_{-\infty}}{n + \langle \x', \s' \rangle} > 0$ completes the
  proof.
\end{proof}

\begin{lem}[$\ell_\infty$ version of shifted scaling {\cite{gonzaga2004properties}}] \label{lem:shiftedscal} 
Suppose $\x, \z \in \Fcal^0_p$ and $\| \X^{- 1} ( \z - \x ) \|_{\infty} \leq \rho \in
  [0, \frac{1}{4}]$. Then
\begin{align}
    & \| ( \X^{- 1} \Z \tmP_{\A\Z} \X^{- 1} \Z - \tmP_{\A\X} )
    \tmv \|
    \leq{}  3 \rho \| \tmP_{\A\X} \tmv \| \nonumber
  \end{align}
 
for all  $\tmv \in
  \mathbb{R}^n$
\end{lem}

  \subsection{Proof of Lemma \ref{lem:boundedprox}}
  Denote $\Delta := \sqrt{1 - \sqrt{1 - 2 \lambda^2}} \leq 1$. 
    For any $\x^1, \x^2 \in \mathcal{P}_{\mu}$, \eqref{ineq:dist-to-central} of \textbf{Lemma \ref{lem:pathproxbnd}} suggests that 
    \begin{align}
    \| \x^1 - \x^{\mu}
    \|_{\x^\mu} \leq \Delta \quad \text{and} \quad 
    \| \x^2 - \x^{\mu} 
    \|_{\x^\mu} \leq \Delta.
  \end{align}
By triangle inequality, the scaled distance between $\x^1$ and $\x^2$ is bounded by
\begin{align*}
    \| \x^1 - \x^2
    \|_{\x^\mu}
    \leq \| \x^1 - \x^\mu
    \|_{\x^\mu} + \| \x^2 - \x^\mu
    \|_{\x^\mu}
    \leq 2 \Delta.
  \end{align*}
  Squaring both sides and partitioning the coordinates by $\mathcal{B},
  \mathcal{N}$,
  \begin{align*}
    \|  \x^1 - \x^2 
    \|^2_{\x^\mu}
    =  \| ( \X_{\mathcal{N}}^{\mu} )^{- 1} ( \x_{\mathcal{N}}^1 - \x_{\mathcal{N}}^2 )
    \|^2 + \| ( \X_{\mathcal{B}}^{\mu} )^{- 1} ( \x_{\mathcal{B}}^1 - \x_{\mathcal{B}}^2 )
    \|^2 
    \leq 4 \Delta^2
  \end{align*}
  where $\X_{\mathcal{N}}^{\mu} = \text{diag}(\x_{\mathcal{N}}^{\mu})$ and $\X_{\mathcal{B}}^{\mu} = \text{diag}(\x_{\mathcal{B}}^{\mu})$. Invoke the lower bound from \tmtextbf{Lemma \ref{lem:cplb}},
\begin{equation}
    	\tfrac{1}{\gamma^2
     \mu^2} \| \x_{\mathcal{N}}^1 - \x_{\mathcal{N}}^2 
    \|^2 \leq \| ( \X_{\mathcal{N}}^{\mu} )^{- 1} ( \x_{\mathcal{N}}^1 - \x_{\mathcal{N}}^2 )
    \|^2  \leq 4 \Delta^2 \leq 4
    \end{equation}
  and thus $\| \x_{\mathcal{N}}^1 - \x_{\mathcal{N}}^2 \| \leq 2
  \gamma \mu$. Next, we bound the term
  $\| \x_{\mathcal{B}}^1 - \x_{\mathcal{B}}^2 \|$.  Since
  \begin{align}
    \A_{\mathcal{B}} \x_{\mathcal{B}}^1 + \A_{\mathcal{N}} \x_{\mathcal{N}}^1
    = \tmb \quad \text{and} \quad 
    \A_{\mathcal{B}} \x_{\mathcal{B}}^2 + \A_{\mathcal{N}} \x_{\mathcal{N}}^2
    = \tmb,
  \end{align}
we subtract them to get $\A_{\mathcal{B}} ( \x_{\mathcal{B}}^1 -
  \x_{\mathcal{B}}^2 ) = - \A_{\mathcal{N}} ( \x_{\mathcal{N}}^1 -
  \x_{\mathcal{N}}^2 )$, which, associated with \ref{A2}, implies that 
  \begin{align*}
    \| \x_{\mathcal{B}}^1 - \x_{\mathcal{B}}^2 \| = \|
    ( \A_{\mathcal{B}}^{\top} \A_{\mathcal{B}} )^{- 1}
    \A_{\mathcal{B}}^{\top} \A_{\mathcal{N}} ( \x_{\mathcal{N}}^1 -
    \x_{\mathcal{N}}^2 ) \|
    \leq \|( \A_{\mathcal{B}}^{\top} \A_{\mathcal{B}} )^{- 1}
    \A_{\mathcal{B}}^{\top}\| \cdot \| \A_{\mathcal{N}}\| \cdot \| \x_{\mathcal{N}}^1 -
    \x_{\mathcal{N}}^2 \|.
  \end{align*}
Using $\|\U\|^2 = \|\U\U^\top\|$,
\[
\|( \A_{\mathcal{B}}^{\top} \A_{\mathcal{B}} )^{- 1}
    \A_{\mathcal{B}}^{\top}\|^2 = \|( \A_{\mathcal{B}}^{\top} \A_{\mathcal{B}} )^{- 1}
    \A_{\mathcal{B}}^{\top}  \A_{\mathcal{B}} (\A_{\mathcal{B}}^{\top} \A_{\mathcal{B}} )^{- 1}\|= 
 \|(\A_{\mathcal{B}}^{\top} \A_{\mathcal{B}})^{- 1}\| = \tfrac{1}{\lambda_{\min} (
     \A_{\mathcal{B}}^{\top} \A_{\mathcal{B}} )}.
\]
Finally, we put things together to get
  \[ \| \x^1 - \x^2 \| \leq \| \x_{\mathcal{N}}^1 -
     \x_{\mathcal{N}}^2 \| + \| \x_{\mathcal{B}}^1 -
     \x_{\mathcal{B}}^2 \| \leq 2 \gamma \mu ( 1 + \tfrac{\|
     \A_{\mathcal{N}} \|}{\sqrt{\lambda_{\min} (
     \A_{\mathcal{B}}^{\top} \A_{\mathcal{B}} )}} ), \]
  and this completes the proof.
  
  \subsection{Proof of Theorem \ref{thm:iterconv}}
  
  Consider $\x^{\mu_k}, \x^{\mu_{k + 1}}$ and let $\s^{\mu_k}, \s^{\mu_{k +
  1}}$ denote the corresponding central path dual solution. We have
  \[ \delta ( \x^{\mu_{k + 1}}, \mu_k ) = \min_{(\y, \s) \in \Fcal_d} \| \tfrac{1}{\mu_k}
     \X^{\mu_{k + 1}} \s - \e \|\leq \| \tfrac{1}{\mu_k}
     \X^{\mu_{k + 1}} \s^{\mu_{k + 1}} - \e \| = \| \tfrac{\mu_{k +
     1}}{\mu_k} \e - \e \| = \tau \sqrt{n} \leq \tfrac{\sqrt{2}}{2} \]
  and $\x_{\mu_{k + 1}} \in \mathcal{P}_{\mu_k}$.
  Since $\delta( \x_k, \mu_k
  ) \leq \frac{\sqrt{2}}{2}$ for all $\x_k$ and $\mu_k$, for all $i \geq k$, the fact that $\x_i \in \mathcal{P}_{\mu_i}$ implies that $\x_i \in \bigcup_{l = k}^\infty \mathcal{P}_{\mu_l}$. Then, for any $i ,j \geq k$,
  $$
  \textstyle \| \x_i - \x_j \| \leq \diam \big( \bigcup_{l = k}^\infty
     \mathcal{P}_{\mu_l} \big) \leq \sum_{l = k}^\infty \diam (\mathcal{P}_{\mu_l})
  $$
  where the last inequality holds because $\x_{\mu_{k + 1}} \in \mathcal{P}_{\mu_k}$ implies that $\mathcal{P}_{\mu_k} \cap \mathcal{P}_{\mu_{k+1}} \neq \varnothing $ holds for all $k$.
 Based on \tmtextbf{Lemma \ref{lem:boundedprox}}, we have
  \begin{align}\label{ineq:cauchy-seq}
      \textstyle \| \x_i - \x_j \| \leq 2 \gamma  ( 1 + \tfrac{\|
     \A_{\mathcal{N}} \|}{\sqrt{\lambda_{\min} (
     \A_{\mathcal{B}}^{\top} \A_{\mathcal{B}} )}} )  \sum_{l = k}^\infty \mu_l
     =2 \gamma  ( 1 + \tfrac{\|
     \A_{\mathcal{N}} \|}{\sqrt{\lambda_{\min} (
     \A_{\mathcal{B}}^{\top} \A_{\mathcal{B}} )}} ) \cdot \mu_k \tau^{- 1} 
  \end{align}
and this completes the proof.

\subsection{Proof of Lemma \ref{lem:euc}}

Using $| 2 x y | \leq \gamma x^2 + \gamma^{- 1} y^2$, for any $0 \leq \gamma
\leq 1$,
\[ (1 - \gamma) \X^2 + (1 - \gamma^{- 1}) \Delta \X^2 \preceq ( \X +
   \Delta \X )^2 \preceq (1 + \gamma) \X^2 + (1 + \gamma^{- 1}) \Delta
   \X^2 \]
Now, since $\A$ is full rank,
\[ (1 - \gamma) \M_{\X} + (1 - \gamma^{- 1}) \A \Delta \X^2 \A^{\top} \preceq
   \M_{\X + \Delta \X} \preceq (1 + \gamma) \M_{\X} + (1 + \gamma^{- 1}) \A
   \Delta \X^2 \A^{\top}. \]
Multiplying
$\M_{\X}^{- 1 / 2}$ on both sides and using $\M_{\X} \succeq \lambda_\x
\cdummy \I$,
\[ \Big[ (1 - \gamma) + (1 - \gamma^{- 1}) \tfrac{\| \A \|^2
   \| \Delta \x \|_{\infty}^2}{\lambda_\x} \Big] \cdummy \I \preceq \M_{\X}^{- 1 / 2} ( \M_{\X +
   \Delta \X} ) \M_{\X}^{- 1 / 2} \preceq \Big[ (1 + \gamma) + (1 +
   \gamma^{- 1}) \tfrac{\| \A \|^2 \| \Delta \x
   \|_{\infty}^2}{\lambda_\x} \Big] \cdummy
   \I. \]
Since $\| \Delta \x \|_{\infty} \leq \tfrac{\sqrt{\lambda_\x}
}{\| \A \|} \beta $, this implies $\tfrac{1}{\lambda_\x }\| \A
\|^2 \| \Delta \x \|^2_{\infty} \leq \beta^2 \leq 1$ and
\[ (1 - \gamma) + (1 - \gamma^{- 1}) \tfrac{\| \A \|^2 \|
   \Delta \x \|_{\infty}^2}{\lambda_\x} = (1 - \gamma^{- 1}) \left( \gamma - \tfrac{\| \A \|^2 \|
   \Delta \x \|_{\infty}^2}{\lambda_\x}\right) \geq
   0 \]
for all $\tfrac{1}{\lambda_{\x} }\| \A \|^2 \| \Delta \x
\|_{\infty}^2 \leq \gamma  \leq 1$. Last taking
$\gamma = \sqrt{\tfrac{1}{\lambda_{\x}}\| \A \|^2 \| \Delta \x
\|_{\infty}^2} \geq \tfrac{1}{\lambda_{\x}}\| \A \|^2 \| \Delta \x
\|_{\infty}^2$, we get
\[ (\gamma - 1)^2 \cdummy \I \preceq \M_{\X}^{- 1 / 2} ( \M_{\X + \Delta
   \X} ) \M_{\X}^{- 1 / 2} \preceq (\gamma + 1)^2 \cdummy \I \]
Notice that $\gamma = \frac{1}{\sqrt{\lambda_{\x}}}
   \| \A \| \| \Delta \x \|_{\infty} \leq \frac{1}{\sqrt{\lambda_{\x}}}
   \| \A \| \| \Delta \x \| \leq \beta \leq 1$, we have
\[
(1 - \beta)^2 \cdummy \I \preceq \M_{\X}^{- 1 / 2} ( \M_{\X + \Delta
   \X} ) \M_{\X}^{- 1 / 2} \preceq (\beta + 1)^2 \cdummy \I
\]
and thus
\[ \kappa ( \M_{\X}^{- 1 / 2} ( \M_{\X + \Delta \X} )
   \M_{\X}^{- 1 / 2} ) \leq ( \tfrac{\beta + 1}{\beta - 1}
   )^2 , \]
which completes the proof.

\subsection{Proof of Theorem \ref{thm:condnum}}

  By \eqref{ineq:cauchy-seq} in the proof of \tmtextbf{Theorem \ref{thm:iterconv}}, for all $j \geq k \geq k^\star$, we have $\mu_k \leq
  \frac{\sqrt{\lambda^{\star}} \tau}{16 \gamma \| \A \|} \cdummy
  ( 1 + \tfrac{\| \A_{\mathcal{N}} \|}{\lambda_{\min} (
  \A_{\mathcal{B}}^{\top} \A_{\mathcal{B}} )} )^{- 1}$ and
  \[ \| \x_j - \x_{k} \| \leq \tfrac{2 \gamma
     \mu_k}{\tau} ( 1 + \tfrac{\| \A_{\mathcal{N}}
     \|}{\sqrt{\lambda_{\min} ( \A_{\mathcal{B}}^{\top}
     \A_{\mathcal{B}} )}} ) \leq \tfrac{\sqrt{\lambda^{\star}}}{8
     \| \A \|}. \]
In particular, $\| \x_k - \x^{\star} \| \leq \tfrac{\sqrt{\lambda^{\star}}}{8
     \| \A \|}$ and using \tmtextbf{Lemma \ref{lem:euc}}, $\M_{\X_k}^{- 1 / 2} \M_{\X^\star}  \M_{\X_k}^{- 1 / 2} \preceq (\beta + 1)^2 \cdummy \I$, we have $\M_{\X_k} \succeq \frac{1}{4}\lambda^\star$.
Finally, invoking \tmtextbf{Lemma \ref{lem:euc}} again with $\beta = \tfrac{1}{2}, \x = \x_k,  \lambda_\x = \frac{1}{4}\lambda^\star$ and $\Delta \x = \x_k - \x_j$ completes the proof.

\subsection{Proof of Theorem \ref{thm:delayedscal}}
\label{app:proof-delay}
Let $\tmv = \X ( \tfrac{1}{\mu} \tmc - \X^{- 1} \e )$, $
  \theta = \| \z - \x \|_{\x, \nu}$,  and define the partition
  $\mathcal{L} \assign \{ j : x_j \geq \nu \}, \mathcal{S} \assign \{ j : x_j
  < \nu \}$. Recall that
  \[ \Delta \x = - \X ( \I - \X \A^{\top} ( \A \X^2 \A^{\top}
     )^{- 1} \A \X ) \tmv . \]
  We show how to find $\Delta \hat{\x}$ combining shifted scaling
  and preconditioning. Let $\w \in \mathbb{R}^n$ be such that $\w_{\mathcal{S}} =
  \x_{\mathcal{S}}$ and $\w_{\mathcal{L}} = \z_{\mathcal{L}}$, and define 
  \[ \Delta \tilde{\x} = - \W ( \I - \W \A^{\top} (
     \A \W^2 \A^{\top} )^{- 1} \A \W ) \W \X^{- 1} \tmv. \]
  By the definition of $\nu$-thresholded scaled distance,
  \[ \| \X^{- 1} ( \x - \w ) \|_{\infty} = \| \X^{-
     1}_{\mathcal{L}} ( \x_{\mathcal{L}} - \z_{\mathcal{L}} )
     \|_{\infty} \leq \| \z - \x \|_{\x, \nu}  \leq \tfrac{1}{4}. \]
Then, we apply \tmtextbf{Lemma \ref{lem:shiftedscal}} and deduce that
  \begin{align}
    & \| \X^{- 1} ( \Delta \x - \Delta \tilde{\x}
    ) \| \nonumber\\
    ={} & \| ( \I - \X \A^{\top} ( \A \X^2 \A^{\top} )^{-
    1} \A \X ) \tmv - \X^{- 1} \W ( \I - \W \A^{\top} ( \A
    \W^2 \A^{\top} )^{- 1} \A \W ) \W \X^{- 1} \tmv \|
    \nonumber\\
    \leq{} & 3 \theta \| ( \I - \X \A^{\top} ( \A \X^2
    \A^{\top} )^{- 1} \A \X ) \tmv \| \nonumber\\
    ={} & 3 \theta \delta ( \x, \mu ), \nonumber
  \end{align}
  where the last equality uses the fact $\delta ( \x, \mu ) =
  \| ( \I - \X \A^{\top} ( \A \X^2 \A^{\top} )^{- 1} \A
  \X ) \X ( \tfrac{1}{\mu} \tmc - \X^{- 1} \e ) \|$ from \textbf{Definition \ref{def:proximity}}. Plugging in condition $\theta \leq \tfrac{\varepsilon}{6 \delta (
  \x, \mu )}$ gives $\| \X^{- 1} ( \Delta \x - \Delta
  \tilde{\x} ) \| \leq \tfrac{\varepsilon}{2}$.  \\
Next consider solving $( \A \W^2 \A^{\top} )^{- 1} \A \W (
  \W \X^{- 1} ) \tmv$. We have $\| \w - \z \| = \| \w_{\mathcal{S}} - \z_{\mathcal{S}}
     \| \leq \theta$ and for $\theta \leq \frac{\sqrt{\lambda_{\z}}}{2 \| \A \|}$, we invoke \tmtextbf{Lemma \ref{lem:euc}} with $\beta = \tfrac{1}{2}$ to get $\A \W^2 \A^{\top}
  \succeq \tfrac{1}{4} \A \Z^2 \A^{\top} \succeq \frac{\lambda_{\z}}{4}$ and
  \[ \kappa ( ( \A \Z^2 \A^{\top} )^{- 1 / 2} \A \W^2
     \A^{\top} ( \A \Z^2 \A^{\top} )^{- 1 / 2} ) \leq 9, \]
  Finally, we analyze the effect of inexactness. Suppose we obtain
  $\hat{\tmd}$ satisfying
  \begin{align}
    \A \W^2 \A^{\top} \hat{\tmd} = \A \W ( \W \X^{- 1} ) \tmv
    +\tmmathbf{\zeta}_{\w} \nonumber
  \end{align}
and the inexact solution $\Delta \hat{\x}$ is chosen as
\begin{align}
    \Delta \hat{\x} ={} & \W^2 \X^{- 1} \tmv - \W^2 \A^{\top}
    \hat{\tmd} \nonumber\\
    ={} & \W^2 \X^{- 1} \tmv - \W^2 \A^{\top} ( \A \W^2 \A^{\top}
    )^{- 1} ( \A \W ( \W \X^{- 1} ) \tmv
    +\tmmathbf{\zeta}_{\w} ), \nonumber
\end{align}
which implies
  \[ \| \X^{- 1} ( \Delta \tilde{\x} - \Delta
     \hat{\x} ) \| = \| \X^{- 1} \W^2 \A^{\top}
     ( \A \W^2 \A^{\top} )^{- 1} \tmmathbf{\zeta}_{\w} \|
     \leq \tfrac{\| \tmmathbf{\zeta}_{\w} \| \cdummy \| \X^{-
     1} \w \|_{\infty}}{\sqrt{\lambda_{\min} ( \A \W^2 \A^{\top}
     )}} \leq \tfrac{2 (\theta + 1)}{\sqrt{\lambda_{\z}}} \|
     \tmmathbf{\zeta}_{\w} \| \]
where the last inequality uses $\| \X^{- 1} \w \|_{\infty} =
  \| \X^{- 1} ( \w - \x ) + \e \|_{\infty} \leq \theta +
  1$. Using an iterative solver, such as conjugate gradient, and preconditioned by $\A \Z^2 \A^{\top}$, we can achieve
  $\| \X^{- 1} ( \Delta \tilde{\x} - \Delta
  \hat{\x} ) \| \leq \frac{\varepsilon}{2}$ in $\mathcal{O} ( \log ( \tfrac{4(\theta + 1)}{\varepsilon \sqrt{\lambda_{\z}}} ) 
     ) = \mathcal{O} ( \log ( \tfrac{1}{\varepsilon } ) 
     )$
iterations. Finally, we put things together to get 
  $$\| \X^{- 1} ( \Delta \x - \Delta
  \hat{\x} ) \| \leq \| \X^{- 1} ( \Delta
  \tilde{\x} - \Delta \hat{\x} ) \| + \|
  \X^{- 1} ( \Delta \x - \Delta \tilde{\x} )
  \| \leq \varepsilon .$$
  This completes the proof.
  
  \begin{rem}
  	We remark that $\Delta \tilde{\x}$ may not satisfy $\A\Delta \tilde{\x} = \0$ due to inexactness and will violate primal feasibility. We address this issue in \textbf{Section \ref{app:inexact}}.
  \end{rem}

\section{Inexact Primal Interior Point Method and Proof of Theorem \ref{thm:informalinexact}} \label{app:inexact}

\subsection{Algorithm Design}
This section provides convergence analysis of \tmtextbf{Algorithm \ref{alg:pipm-precond}}.
We start by defining inexactness. Note that inexact direction $\Delta
\x$ can be re-written as
\begin{align}
  \Delta \x ={} & - \X ( \I - \X \A^{\top} ( \A \X^2 \A^{\top}
  )^{- 1} \A \X ) ( \tfrac{1}{\mu} \X \tmc - \e )
  \nonumber\\
  ={} & - ( \tfrac{1}{\mu} \X^2 \tmc - \x ) + \X^2 \A^{\top} (
  \A \X^2 \A^{\top} )^{- 1} ( \tfrac{1}{\mu} \A \X^2 \tmc - \A \x
  ) \nonumber\\
  ={} & - ( \tfrac{1}{\mu} \X^2 \tmc - \x ) + \X^2 \A^{\top} \tmd,
  \nonumber
\end{align}

where the normal equation $\A \X^2 \A^{\top} \tmd = \tfrac{1}{\mu} \A \X^2
\tmc - \A \x$ is solved inexactly. Suppose an error $\tmmathbf{\zeta}_{\x}$ is incurred:
\begin{equation} \label{eqn:inexact-neq}
	\A \X^2 \A^{\top} \hat{\tmd} = \tfrac{1}{\mu} \A \X^2 \tmc - \A \x
   +\tmmathbf{\zeta}_{\x}
\end{equation}
If we directly substitute back, $\Delta \hat{\x} = - (
\tfrac{1}{\mu} \X^2 \tmc - \x ) + \X^2 \A^{\top} \hat{\tmd}$ and
\begin{align}
  \A \Delta \hat{\x} ={} & - \A ( \tfrac{1}{\mu} \X^2 \tmc - \x
  ) + \A \X^2 \A^{\top} \hat{\tmd} \nonumber\\
  ={} & - \A ( \tfrac{1}{\mu} \X^2 \tmc - \x ) + \tfrac{1}{\mu} \A
  \X^2 \tmc - \A \x +\tmmathbf{\zeta}_{\x} ={}\tmmathbf{\zeta}_{\x} . \nonumber
\end{align}

In other words, the error in solving the normal equation propagates and results in
primal infeasibility, which is undesirable. There is a simple fix: take
arbitrary basis matrix $\B$ of $\A$ and denote the column partition of $\A$ as $(B, N)$. Then we can take $\tmmathbf{\lambda}_{\x}$ such that
$\tmmathbf{\lambda}_{\x, B} = \B^{- 1} \tmmathbf{\zeta}_{\x}$ and
$\tmmathbf{\lambda}_{\x, N} = \0$. Taking $\Delta \hat{\x} = -
( \tfrac{1}{\mu} \X^2 \tmc - \x ) + \X^2 \A^{\top} \hat{\tmd}
+\tmmathbf{\lambda}_{\x}$, this gives $\A \Delta \hat{\x}= \0$. For
brevity of analysis, we take
\[ \tmmathbf{\lambda}_{\x} \assign \argmin_{\A \tmmathbf{\lambda}=
   -\tmmathbf{\zeta}_{\x}}  \| \tmmathbf{\lambda} \| = - \A^{\top} ( \A
   \A^{\top} )^{- 1} \tmmathbf{\zeta}_{\x} . \]
From now on we formally define $\Delta \hat{\x} \assign - (
\tfrac{1}{\mu} \X^2 \tmc - \x ) + \X^2 \A^{\top} \hat{\tmd}
+\tmmathbf{\lambda}_{\x}$.\\

To analyze the convergence of \textbf{Algorithm \ref{alg:pipm-precond}}, we first present auxiliary results in \textbf{Section \ref{app:inexact-auxi}}. We prove the complexity of \tmtextbf{Algorithm \ref{alg:pipm-precond}} in
\tmtextbf{Section \ref{app:inexact-analysis}}. Without loss of generality, in this section, we drop the iteration index $k$ and only consider two consecutive IPM iterations $\x, \x^+$ with barrier parameters $\mu, \mu^+$ respectively. Given $\x \in \mathcal{F}^0_p$, we define the projected dual solution
\[ \s ( \x, \mu ) \assign \argmin_{( \y, \s ) \in
   \mathcal{F}_d} ~ \| \tfrac{1}{\mu} \X \s - \e \| \]
and let $\z \assign \tfrac{1}{\mu} \X \s ( \x, \mu )$.

\subsection{Auxiliary Results}\label{app:inexact-auxi}

\begin{lem}[Primal step \cite{roos1992polynomial}] 
\label{lem:primal-step}
The following identities hold
  \begin{align}
    \tmP_{\A\X} ( \tfrac{1}{\mu} \X \tmc - \e ) ={}  & \tfrac{1}{\mu} \X
    \s ( \x, \mu ) - \e \label{eqn:primal-step-1}\\
    \Delta \x ={} & \x - ( \tfrac{1}{\mu} \X^2 \s ( \x, \mu
    ) - \x ) = 2\x - \X \z \label{eqn:primal-step-2}
  \end{align}
\end{lem}

\begin{lem}[Path-following \cite{roos1992polynomial}] \label{lem:path-following}
If $\mu^+ = (1 - \tau) \mu, \tau \in (0, 1)$,
  then $ \delta ( \x, \mu^+ ) \leq \tfrac{1}{1 - \tau} ( \delta
     ( \x, \mu ) + \tau \sqrt{n} )$.
\end{lem}

\subsection{Inexact Analysis} \label{app:inexact-analysis}
\begin{lem} \label{lem:inexactstep}
  Suppose $\delta ( \x, \mu ) \leq \frac{1}{2}$, we have
  \[ \x^+ = \X ( 2 \e - \z +\tmmathbf{\psi}_{\x} ) \quad \text{and}
     \quad \delta ( \x^+, \mu ) \leq \sqrt{2} \delta ( \x, \mu
     )^2 + \sqrt{2 n} \| \tmmathbf{\psi}_{\x} \|, \]
  where $\tmmathbf{\psi}_{\x} \assign \X \A^{\top} ( \A \X^2 \A^{\top}
  )^{- 1} \tmmathbf{\zeta}_{\x} + \X^{- 1} \tmmathbf{\lambda}_{\x}$.
\end{lem}

\begin{proof}
  By definition of the inexact primal IPM step, we successively deduce that
  \begin{align}
    \x^+ ={} & \x + \Delta \hat{\x}\nonumber\\
    ={} & \x + \Delta \x + \X^2 \A^{\top} ( \A \X^2 \A^{\top}
    )^{- 1} \tmmathbf{\zeta}_{\x} +\tmmathbf{\lambda}_{\x} \nonumber\\
    ={} & 2 \x - \X \z + \X^2 \A^{\top} ( \A \X^2 \A^{\top} )^{- 1}
    \tmmathbf{\zeta}_{\x} +\tmmathbf{\lambda}_{\x} \label{eqn:proof-lem-c-3-1} \\
    ={} & \X ( 2 \e - \z + \X \A^{\top} ( \A \X^2 \A^{\top} )^{-
    1} \tmmathbf{\zeta}_{\x} + \X^{- 1} \tmmathbf{\lambda}_{\x} )
    \nonumber\\
    ={} & \X ( 2 \e - \z +\tmmathbf{\psi}_{\x} ), \nonumber
  \end{align}
where \eqref{eqn:proof-lem-c-3-1} applies \eqref{eqn:primal-step-2}. In view of the proximity measure, we have
  \begin{align}
    \delta ( \x^+, \mu ) =  \| \tfrac{1}{\mu} \X^+ \s (
    \x^+, \mu ) - \e \| 
    \leq{} & \| \tfrac{1}{\mu} \X^+ \s ( \x, \mu ) - \e \|
    ={} \| \X^+ \X^{- 1} \z - \e \|. \nonumber
  \end{align}
 
  Plugging in $\x^+ = \X ( 2 \e - \z +\tmmathbf{\psi}_{\x} )$,
  \[ \X^+ \X^{- 1} \z = \X^{- 1} \X ( 2 \I - \Z +\tmmathbf{\Psi}_{\x}
     ) \z = 2 \z - \Z \z +\tmmathbf{\Psi}_{\x} \z, \]
  where $\tmmathbf{\Psi}_{\x} := \diag( \tmmathbf{\psi}_{\x}
  )$. Hence, we deduce that 
  \begin{align}
    \delta ( \x^+, \mu )^2 ={} & \| \X^+ \X^{- 1} \z - \e
    \|^2 \nonumber\\
    ={} & \| 2 \z - \Z \z - \e +\tmmathbf{\Psi}_{\x} \z \|^2
    \nonumber\\
    ={} & \| 2 \z - \Z \z - \e \|^2 + \langle
    2\tmmathbf{\Psi}_{\x} \z, 2 \z - \Z \z - \e \rangle + \|
    \tmmathbf{\Psi}_{\x} \z \|^2 \nonumber\\
    \leq{} & 2 \delta ( \x, \mu )^4 + 2
    \| \tmmathbf{\Psi}_{\x} \z \|^2,\label{eqn:proof:lemc31-1}
  \end{align}
  where \eqref{eqn:proof:lemc31-1} uses Cauchy's inequality.   Taking square root on both sides and using $\sqrt{x + y} \leq \sqrt{x} + \sqrt{y}$, 
  \[ \delta ( \x^+, \mu ) \leq \sqrt{2} \delta ( \x^+,
     \mu )^2 + \sqrt{2} \| \tmmathbf{\Psi}_{\x} \z
     \| \leq \sqrt{2} \delta ( \x^+, \mu )^2 + (\sqrt{2n} + 1) \| \tmmathbf{\psi}_{\x} \|, \]
  where we use the relation $\| \z - \e \| \leq \frac{1}{2}$ and $\| \z
  \| \leq \|\e\| + \frac{1}{2} = \sqrt{n} +\frac{1}{2}$.  This completes the proof.
\end{proof}

\begin{lem} \label{lem:inexacterr}
  Suppose $\delta ( \x, \mu ) \leq \tfrac{1}{2}$, then for all $\mu
  \leq \frac{\sqrt{\lambda^{\star}} \tau}{16 \gamma} \cdummy (
  \tfrac{\| \A_{\mathcal{N}} \|}{\sqrt{\lambda_{\min} (
  \A_{\mathcal{B}}^{\top} \A_{\mathcal{B}} )}} + 1 )^{- 1}$  
  \[ \| \tmmathbf{\psi}_{\x} \| \leq (\tfrac{1}{ \alpha \mu \sigma_{\min}
   ( \A )}+\tfrac{2}{\sqrt{\lambda^\star}})  \tfrac{1}{ \sigma_{\min}
   ( \A )}  \| \tmmathbf{\zeta}_{\x} \|. \]
\end{lem}

\begin{proof}
Recall that $\tmmathbf{\psi}_{\x} = \X \A^{\top} ( \A \X^2 \A^{\top}
)^{- 1} \tmmathbf{\zeta}_{\x} + \X^{- 1} \tmmathbf{\lambda}_{\x}$ and we
first consider $\tmmathbf{\lambda}_{\x}$. First note that
\[ \| \tmmathbf{\lambda}_{\x} \| = \sqrt{\langle
   \tmmathbf{\zeta}_{\x}, ( \A \A^{\top} )^{- 1}
   \tmmathbf{\zeta}_{\x} \rangle} \leq \tfrac{1}{\sqrt{\lambda_{\min}
   ( \A \A^{\top} )}} \| \tmmathbf{\zeta}_{\x} \| =
   \tfrac{1}{\sigma_{\min} ( \A )} \| \tmmathbf{\zeta}_{\x}
   \| \]
Then using $\delta ( \x, \mu ) \leq \tfrac{1}{2}$ and
\tmtextbf{Lemma \ref{lem:app:cplb}}, $\| \x \|_{- \infty} \geq \alpha \mu$ and
$\| \X^{- 1} \tmmathbf{\lambda}_{\x} \| \leq \tfrac{1}{\alpha \mu
\sigma_{\min} ( \A )} \| \tmmathbf{\zeta}_{\x} \|$. Next
we consider $\X \A^{\top} ( \A \X^2 \A^{\top} )^{- 1}
\tmmathbf{\zeta}_{\x}$. Similarly, we have
\[ \| \X \A^{\top} ( \A \X^2 \A^{\top} )^{- 1}
   \tmmathbf{\zeta}_{\x} \| = \sqrt{\langle \tmmathbf{\zeta}_{\x},
   ( \A \X^2 \A^{\top} )^{- 1} \tmmathbf{\zeta}_{\x} \rangle}
   \leq \tfrac{1}{\sqrt{\lambda_{\min} ( \A \X^2 \A^{\top} )}}
   \| \tmmathbf{\zeta}_{\x} \| \leq \tfrac{2}{\sqrt{\lambda^\star}} \| \tmmathbf{\zeta}_{\x} \|, \]
where we again invoked \textbf{Lemma \ref{lem:euc}} with $\beta = \frac{1}{2}$ to obtain $\lambda_{\min} ( \A \X^2 \A^{\top} ) \succeq \frac{\lambda^\star}{4} \I$ . Putting things together,
\[ \| \tmmathbf{\psi}_{\x} \| = \| \X \A^{\top} ( \A \X^2
   \A^{\top} )^{- 1} \tmmathbf{\zeta}_{\x} + \X^{- 1}
   \tmmathbf{\lambda}_{\x} \| \leq  \big(\tfrac{1}{ \alpha \mu \sigma_{\min}
   ( \A )}+\tfrac{2}{\sqrt{\lambda^\star}}\big)  \| \tmmathbf{\zeta}_{\x} \| \]
and this completes the proof.	
\end{proof}

\subsection{Proof of Theorem \ref{thm:informalinexact}}

\begin{thm}[\textbf{Theorem \ref{thm:informalinexact}}]
  Under \ref{A1} and \ref{A2} and assume $n > 100$, suppose primal IPM starts from
  \[ \mu_1 < \tfrac{\sqrt{\lambda^{\star}}}{40 \sqrt{n} \gamma} (
     \tfrac{\| \A_{\mathcal{N}} \|}{\sqrt{\lambda_{\min} (
     \A_{\mathcal{B}}^{\top} \A_{\mathcal{B}}} )} + 1 )^{- 1}
  \]
  such that $\delta ( \x_1, \mu_1 ) \leq \frac{1}{2}$ and suppose $\tau = \frac{1}{10 \sqrt{n}}$. Then  $\delta ( \x_k, \mu_k ) \leq
  \tfrac{1}{2}$ for all $k$ and the arithmetic complexity of finding an $\varepsilon$-optimal primal solution is  $\mathcal{O} ( n^{2.5} \log^2 ( \tfrac{1}{\varepsilon}
     ) )$.
\end{thm}

\begin{proof}
 Putting all the results together,
  \begin{align}
    \x^+ ={} & \X ( 2 \e - \z +\tmmathbf{\psi}_{\x} ) \label{eqn:final-1}\\
    \| \tmmathbf{\psi}_{\x} \| \leq{} & (\tfrac{1}{ \alpha \mu \sigma_{\min}
   ( \A )}+\tfrac{2}{\sqrt{\lambda^\star}})   \|
    \tmmathbf{\zeta}_{\x} \| \label{eqn:final-2}\\
    \delta ( \x^+, \mu ) \leq{} & \sqrt{2} \delta ( \x, \mu
    )^2 + (\sqrt{2 n}  + 1)\| \tmmathbf{\psi}_{\x} \| \label{eqn:final-3}\\
    \delta ( \x, \mu^+ ) \leq{} & \tfrac{1}{1 - \tau} ( \delta
    ( \x, \mu ) + \tau \sqrt{n} )  \label{eqn:final-4}
  \end{align}
  Assume we let \[ \| \tmmathbf{\zeta}_{\x} \| \leq \tfrac{\alpha \mu \sigma_{\min} 
   ( \A )\sqrt{\lambda^\star}}{40\sqrt{3n}(\sqrt{\lambda^\star}+2\alpha \mu \sigma_{\min}
   ( \A ))} \]
and by \eqref{eqn:final-2}  $\| \tmmathbf{\psi}_{\x} \| \leq
  \frac{1}{40 \sqrt{3 n}}$. With \textbf{Theorem \ref{thm:condnum}} it takes  $\mathcal{O} ( \log \tfrac{1}{\mu} )$ iterations for an iterative solver to converge.
  Since $\x, \mu$ satisfies $\delta ( \x, \mu ) \leq
  \tfrac{1}{2}$. We have $\| \z - \e \|_{\infty} \leq \| \z -
  \e \| \leq \tfrac{1}{2}$ and from \eqref{eqn:final-1} we deduce that
  \[ 2 \e - \z +\tmmathbf{\psi}_{\x} \geq ( \tfrac{1}{2} - \|
     \tmmathbf{\psi}_{\x} \| ) \e \geq 0.1 \e. \]
     
Therefore $\x^+ \in \mathcal{F}_p^0$. Next fix $\tau = 1 / ( 10 \sqrt{n}
  )$ and for $n \geq 100$
  \begin{align}
    \delta ( \x^+, \mu^+ ) \leq{} & \tfrac{1}{1 - \tau} ( \delta
    ( \x^+, \mu ) + \tau \sqrt{n} ) \nonumber\\
    \leq{} & \tfrac{1}{1 - \tau} ( \sqrt{2} \delta ( \x, \mu )^2
    + (\sqrt{2 n } + 1) \| \tmmathbf{\psi}_{\x} \| ) + \tfrac{\tau
    \sqrt{n}}{1 - \tau} \nonumber\\
    ={} & \tfrac{\sqrt{2}}{1 - \tau} \delta ( \x, \mu )^2 +
    \tfrac{\sqrt{2 n} + 1}{1 - \tau} \| \tmmathbf{\psi}_{\x} \| +
    \tfrac{\tau \sqrt{n}}{1 - \tau} \nonumber\\
    \leq{} & \tfrac{3}{2} \delta ( \x, \mu )^2 + \sqrt{3 n} \|
    \tmmathbf{\psi}_{\x} \| + \tfrac{1}{10} \nonumber
  \end{align}
Now that $\delta ( \x, \mu ) \leq \frac{1}{2}$, we have
  \begin{align}
    \delta ( \x^+, \mu^+ ) \leq{} & \tfrac{3}{8} + \sqrt{3 n} \|
    \tmmathbf{\psi}_{\x} \| + \tfrac{1}{10} 
    \leq{}  \tfrac{19}{40} + \sqrt{3 n} \| \tmmathbf{\psi}_{\x} \| \leq \tfrac{19}{40} + \tfrac{1}{40} = \tfrac{1}{2}.
    \nonumber
  \end{align}
According to  \textbf{Theorem \ref{thm:iterconv}}, we have $\|\x_k - \x^\star \| = \mathcal{O}(\mu_k) = \mathcal{O}((1 - \tau)^k)$, and the total outer iterations are 
$$
\mathcal{O}(\log( 1- \tau) \cdot \log \varepsilon)  = \mathcal{O}(\sqrt{n} \log(\tfrac{1}{\varepsilon}) ).
$$
Besides, in each iteration, the preconditioned iterative algorithm needs to take $\log(\frac{n}{\varepsilon})$ iterations to find an inexact solution with accuracy $\varepsilon$, and the cost of each iteration is the matrix-vector product and is thus $n^2$. This completes the proof. 
\end{proof}



\end{document}